\documentclass[10pt]{amsart}

\usepackage{amsmath}
\usepackage{amssymb}
\usepackage{times}
\usepackage{mathrsfs}
\usepackage{graphicx}
\usepackage{placeins}
\usepackage{subfigure}
\usepackage{wrapfig}

\begin{document}
\numberwithin{equation}{section}

\def\1#1{\overline{#1}}
\def\2#1{\widetilde{#1}}
\def\3#1{\widehat{#1}}
\def\4#1{\mathbb{#1}}
\def\5#1{\frak{#1}}
\def\6#1{{\mathcal{#1}}}

\newcommand{\bC}{{\mathbb{C}}}
\newcommand{\bP}{{\mathbb{P}}}
\newcommand{\bF}{{\mathbb{F}}}
\newcommand{\bZ}{{\mathbb{Z}}}
\newcommand{\bR}{{\mathbb{R}}}
\newcommand{\si}{{\sigma}}
\newcommand{\sid}{{\check{\sigma}}}
\newcommand{\Si}{{\Sigma}}
\newcommand{\D}{{\Delta}}
\newcommand{\al}{{\alpha}}
\newcommand{\la}{{\lambda}}
\newcommand{\G}{{\Gamma}}
\newcommand{\w}{{\widetilde{\D}}}
\newcommand{\tX}{{\widetilde{X}}}
\newcommand{\tp}{{\widetilde{p}}}
\newcommand{\tq}{{\widetilde{q}}}
\newcommand{\tS}{{\widetilde{\Si}}}
\newcommand{\f}{{\varphi}}
\newcommand{\li}{{\mathscr{L}}}
\newcommand{\lp}{{\mathcal{L}_{n,d}(\sum_{i=0}^{k}m_ip_i)}}
\newcommand{\lb}{{\mathcal{L}_{n,d}(2^{k+1})}}
\newcommand{\de}{{\mathcal{D}}}
\newcommand{\s}{{\mathcal{S}}}
\newcommand{\pP}{{\mathcal{P}}}

\newtheorem{teo}{Theorem}[section]
\newtheorem{lemma}[teo]{Lemma}
\newtheorem{defi}[teo]{Definition}
\newtheorem{prop}[teo]{Proposition}
\newtheorem{cor}[teo]{Corollary}
\newtheorem{pb}[teo]{Problem}

\newtheorem{exer}{Exercise}[section]
\newtheorem{example}{Example}[section]
\newtheorem{remark}[teo]{Remark}

\title{A combinatorial approach to Alexander--Hirschowitz's Theorem based on toric degenerations}

\noindent\author[S. Brannetti]{Silvia Brannetti} \address{S. Brannetti: Dipartimento Di Matematica\\ Universit\`{a} di Roma 3\ \\ Largo San Leonardo Murialdo, 1, 00146 Roma, Italy} \email{brannett@mat.uniroma3.it}

\begin{abstract} We present an alternative proof of the
Alexander--Hirschowitz's Theorem in dimension $3$ using degenerations
of toric varieties.
\end{abstract}

\maketitle

\section{Introduction}
\noindent Secant varieties are classical objects of study in algebraic geometry: given a closed variety $X$ in some projective space $\bP^n$, and given a natural number $k$, the $k$-th secant variety of X is the Zariski closure of the union of all subspaces of $\bP^n$ that are spanned by $k+1$ independent points on $X$. Not enough is known about these varieties, starting from their dimensions. Cases of interest are for instance the secant varieties of Segre embeddings of products of projective spaces, Pl\"{u}cker embeddings of Grassmannians, and Veronese embeddings of projective spaces. The Alexander--Hirschowitz's Theorem (\ref{teorema}) provides a complete answer to the case of all Veronese embeddings of $\bP^n$. We recall that the problem of the secant varieties is equivalent to another classical one in Algebraic Geometry, Polynomial Interpolation (see Remark \ref{interludio}), and in this setting Alexander--Hirschowitz's Theorem solves exhaustively the problem on interpolation with double points in any dimension; its original proof is due to Alexander and Hirschowitz in the nineties, and in 2002 it has been simplified by Chandler \cite{bib:cha}. More recently the proof has been improved by Brambilla and Ottaviani in \cite{bib:OttBra}. Moreover a proof of this theorem has been given in dimension $2$ by Draisma (\cite{bib:draisma}) using a new approach based on tropical geometry, and again in dimension $2$ a different proof of the theorem is by Ciliberto, Miranda, Dumitrescu (\cite{bib:CC}) involving  degenerations of toric varieties, that allow to translate the classical problem into an easy combinatorial one. Indeed tropical geometry is closely related to degenerations of toric varieties, so that the two methods are in fact very much connected. For related concepts see (\cite{bib:SS}). In this paper we want to prove Alexander--Hirschowitz's Theorem in dimension $3$ applying toric degenerations, in line with the approach of \cite{bib:CC}, and exploiting the combinatorial nature of toric varieties.
\section*{Acknowledgements}
\noindent I wish to thank Prof. Ciro Ciliberto for several useful discussions, precious suggestions and patient encouragements. I also thank the referee of the paper for his very detailed report and his advices.

\section{Secant varieties and the Alexander--Hirschowitz's theorem}\label{sez1}

\begin{defi}\label{secant}
Let $X\subset\bP^N$ be an irreducible, non--degenerate projective
variety of dimension $n$ and let $k$ be a positive integer. Take
$k+1$ independent points $p_0,\ldots ,p_k$ of $X$. The span $\langle
p_0,\ldots ,p_k\rangle$ is a subspace of $\bP^N$ of dimension $k$
which is called a $(k+1)$--secant $\bP^k$ of $X$. By {\rm
$Sec_k(X)$} we denote the closure of the union of all
$(k+1)$--secant $\bP^k$'s of $X$. This is an irreducible algebraic
variety which is called the {\rm $k$-th secant variety of $X$}.
\end{defi}

\noindent The study of secant varieties in particular concerns their
dimension, which in most cases is unknown. It is easy to see that
there is a natural upper bound on the dimension:
\begin{equation}\label{dim}
dim(Sec_k(X))\leq \min\{(n+1)(k+1)-1,N\}.
\end{equation}
The right hand side of (\ref{dim}) is the \textit{expected
dimension} of $Sec_k(X)$.

\begin{defi}\label{defect}
X is said to be {\rm $k$--defective} if the strict inequality holds
in (\ref{dim}), and its k--defect is {\rm
$\delta_k=\min\{(n+1)(k+1)-1,N\}-dim(Sec_k(X))$}, otherwise it is
said to be {\rm not $k$--defective}.
\end{defi}
\noindent Our aim in this paper is to prove Alexander--Hirschowitz's
Theorem in dimension $3$, so we are going to recall its terms. Let $V_{n,d}$
be the Veronese variety in $\bP^{{n+d \choose d}-1}$:
\begin{teo}{\bf (Alexander--Hirschowitz)}\label{teorema}
The Veronese variety $V_{n,d}$ is always not $k$--defective, except
in the cases: $$\begin{array}{cccccc}

                   n & any & 2 & 3 & 4 & 4\\

                   d & 2 & 4 & 4 & 4 & 3\\

                   k & 2,\ldots ,n & 5 & 9 & 14 & 7.\\

\end{array}
$$
\end{teo}

\noindent Terracini's lemma (see \cite{bib:interp}) allows us to reduce the problem of
computing the dimension of $Sec_k(X)$ to the study of the dimension
of the span of the linear spaces $T_{X,p_i}, i=0,\ldots ,k$ which
are the tangent spaces to $X$ at $p_0,\ldots ,p_k$, for $p_0,\ldots ,p_k$ general points on $X$. Let us observe
that if the $T_{X,p_i}$ with
$i=0,\ldots ,k$ are independent in $\bP^N$, then their span has maximal
dimension, equal to $(n+1)(k+1)-1$, therefore $X$ is not
$k$--defective. Since our aim is to show
Alexander--Hirschowitz's Theorem in dimension $3$, we are
interested in studying the defectivity of the secant variety of the
Veronese threefold $V_{3,d}$, which is isomorphic to $\bP^3$ via the
Veronese embedding. We want to show that $V_{3,d}$ is not
$k$--defective, for any $k\in\mathbb{N}$ and $d\geq 5$. The cases $d\leq 4$ are well known. In general,
given a variety $X$ as in Definition \ref{secant}, let $\de$ be a
degeneration of $X$ to a union of varieties $X_1,\ldots ,X_m$; fix
$p_0,\ldots ,p_k$ general points on $X$. If we let the $k+1$ points
degenerate in such a way that each point goes to a general point of
some component of the central fiber, we can consider the tangent
spaces to the central fiber at the limit points. We denote by
$T_{X_i,p}$ the tangent space to the component $X_i$ of the central
fiber at $p$, where $p$ is a limit point belonging to $X_i$. If the tangent spaces to the central fiber at the limit points are independent, so are the tangent spaces $T_{X,p_i}, i=0,\ldots ,k$, hence $X$ is not $k$--defective.

\begin{remark}\label{interludio}
{\rm We explain briefly the connections of the problem of secant varieties with
polynomial interpolation. The latter can be stated as follows: given a set of
points $p_0,\ldots ,p_k$ in $\bP^n$,
and assigned relative multiplicities $m_0,\ldots ,m_k\in
\mathbb{N}$, we study the dimension of the linear system $\lp$ of
the hypersurfaces of degree $d$ in $\bP^n$ having at least
multiplicity $m_i$ at $p_i$ for each $i=0,\ldots ,k$. The system is
said to be \textit{special} if it doesn't have the expected dimension, otherwise it
is said to be \textit{non--special}. The
study of interpolation problem consists in classifying the special
systems. This has been completely done for $n=1$, $n=2$ and for any
$n$ when the multiplicities are $m_i=2$ $\forall i$. The latter case
is examined in Alexander--Hirschowitz's Theorem, as noted above. By Terracini's lemma we can say that $Sec_k(V_{n,d})$ has the expected dimension if and only if the system $\lb$ is non--special, indeed the Veronese morphism $v_{n,d}$ allows us to embed
$\bP^n$ in $\bP^{{n+d \choose d}-1}$, and to translate the system of the
hypersurfaces of degree $d$ with $k+1$ double points, to the system of
hyperplanes in $\bP^{{n+d \choose d}-1}$ which are tangent to $V_{n,d}$ at $k+1$ fixed general points. In \cite{bib:CC} the authors prove Alexander--Hirschowitz's Theorem in dimension $2$, using suitable degenerations of the Veronese surface $V_{2,d}$, and reducing the problem to an easy combinatorial one. Even if one could use the same technique in dimension three, we will give a similar
approach that emphasizes the perspective of toric varieties as described below.}

\end{remark}

\subsection{Degenerations of toric varieties}
In this paper we will use degenerations of \textit{toric varieties}. Let us briefly recall how a degeneration of a toric
variety can be described. The interested reader is referred to
\cite{bib:shengda} and \cite{bib:CC} for details. A rational polytope
$\pP$ in $\bR^n$ defines a projective toric variety $X_\pP$ endowed with an
ample line bundle; a subdivision $\G$ of $\pP$ is a partition of $\pP$ into smaller polytopes, i.e. there exist polytopes $\pP_1,\ldots,\pP_l$ of the same dimension as $\pP$, such that $\pP_i\cap\pP_j$ is a face of both (it can be the empty face), and $\bigcup_{i=1}^l\pP_i=\pP$. Then $\Gamma$ is said to be a \textit{regular
subdivision} of $\pP$ if there exists an integral function $F$ defined on $\pP$, which is piecewise linear on the subpolytopes of $\G$ and strictly convex on
$\Gamma$, i.e. for any points $p$ and $q$ in different subpolytopes
of $\Gamma$, $F(tp+(1-t)q)<tF(p)+(1-t)F(q)$ for all $t\in[0,1]$. If $\G$ is regular,
then there exists a degeneration of $X_\pP$ such that the central
fiber is the divisor $X_0$ whose components are the toric varieties
defined by the subpolytopes of $\Gamma$, whereas the total space of
the degeneration is the toric variety associated to the unbounded
polytope $\widetilde{\pP}=\{(x,z)\in\pP\times\bR:z\geq F(x)\}$. In
this construction $X_\pP$ is isomorphic to the general fiber, thus
it degenerates to $X_0$. If $X_i$ is a component of the central fiber endowed with a very ample line bundle, then the corresponding subpolytope $\pP_i$  provides an embedding of $X_i$ in the projective space $\bP^{l_i-1}$, where $l_i$ is the number of integral points in $\pP_i$. Note that if two subpolytopes $\pP_i,\pP_j$ in a subdivision $\G$ of $\pP$ are disjoint, then the varieties they define are embedded in independent linear subspaces of $\bP^{M-1}$, where $M$ is the number of integral points of $\pP$.

\begin{prop}\label{prop}
Let $X\subset\bP^N$ be an irreducible, toric projective
variety of dimension $n$, let $\de$ be a toric degeneration of $X$ to
$X_0=X_1\cup \ldots\cup X_m$, and let $\pP_1,\ldots,\pP_m$ be the corresponding polytopes. If
\begin{itemize}
\item[(i)]there exist indices $j_1,\ldots ,j_h\in\{ 1,\ldots ,m\}$ such that $X_{j_i}$ is not $k_{j_i}$--defective,
\item[(ii)]$\pP_{j_i}\cap \pP_{j_l}=\emptyset$ for any $i\neq l,$ with $i,l=1,\ldots,h$,
    \end{itemize}
then, setting $k+1=\sum_{i=1}^h (k_{j_i}+1)$, $X$ is not $k$--defective.
\end{prop}

\noindent\textit{Proof}. We can assume $j_1=1,\ldots,j_h=h$, so $X_i$ is not $k_i$--defective, for any $i=1,\ldots,h$; fix $p_j^i$, $j=0,\ldots ,k_i$ general points on $X_i$. Then $dim(\langle \bigcup_{j=0,\ldots ,k_i}T_{X_i,p^i_j}\rangle)=(n+1)(k_i+1)-1$. Since the polytopes $\pP_1,\ldots,\pP_h$ are disjoint, then we can conclude that $$dim(\langle \bigcup_{j=0,\ldots ,k_i \atop i=1,\ldots,h}T_{X_i,p^i_j}\rangle)=(n+1)(k+1)-1.$$
We end the proof by observing that the points $p_j^i$ are limits of general points on $X$.

\rightline{$\Box$}

\bigskip

\noindent Now let $X$ be the Veronese threefold $V_{3,d}$ in
$\bP^{N_d}$, where $N_d=\left( {d+3 \atop 3}\right)-1$. We want to prove Alexander--Hirschowitz's Theorem for any $d\geq 5$ and any
$k\in\mathbb{N}$. Let $n_d+1=\left\lfloor \frac{(d+1)(d+2)(d+3)}{24}
\right\rfloor$; we have the following results, whose easy proofs can be left to the reader:

\begin{lemma}\label{expdim}
The expected dimension of $Sec_{n_d}(V_{3,d})$ is\begin{equation}\label{secondagraffa}
\begin{cases}
N_d & \text{ if } d\text{ is odd or }d=6+8k,\text{ for }k\geq 0\\
N_d-1 & \text{ if } d=8k,\text{ for }k\geq 1\\
N_d-2 & \text{ if } d=10+8k,\text{ for }k\geq 0\\
N_d-3 & \text{ if } d=12+8k,\text{ for }k\geq 0.
\end{cases}
\end{equation}
\end{lemma}

\begin{lemma}\label{reduce}
If $V_{3,d}$ is not $n_d$--defective for any $d\geq 5$,
then it is not $k$--defective for any $k<n_d$.
\end{lemma}

\begin{remark}\label{reduceremark}
{\rm Suppose we have proved that $V_{3,d}$ is not $n_d$--defective. If $k>n_d$, in order to see that $V_{3,d}$
is not $k$--defective we have to prove that the dimension of
the span of the tangent spaces at $k+1$ general points is maximal.
By Lemma \ref{expdim} this is obvious when either $d$ is odd,
or when $d\equiv 6 \mod 8$; in the other cases we
have that $Sec_{n_d}(V_{3,d})$ has codimension $1,2$ or $3$
depending on $d\equiv 0,2$ or $4 \mod 8$. If we impose one further
general point in such a way that its tangent space and the tangent
spaces to the points that we have already imposed span the whole
ambient space $\bP^{N_d}$, the dimension of $Sec_{n_d}(V_{3,d})$ is maximal.}
\end{remark}

\begin{figure}[h] \centering
\includegraphics[height=2.8cm]{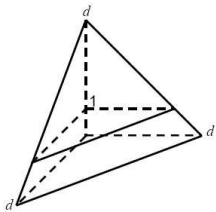}
\caption{\small{Section of $\D_d$ into $\D_{d-1}$ plus one layer of height $1$.}}\label{tetrasezione}
\end{figure}

\begin{figure}[h] \centering
\includegraphics[height=1.4cm]{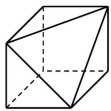}
\caption{\small{The block $\Sigma$.}}\label{blocco}
\end{figure}

\noindent The degeneration of the Veronese threefold we will deal with, is described in the following:
\begin{lemma}\label{degeneration1}
The Veronese threefold $V_{3,d}$ degenerates to $\left({d \atop 3}\right)$ copies of $\bP^1\times\bP^1\times\bP^1$, plus $\frac{d(d+1)}2$ copies of
$\bP^3$, and $\frac{d(d-1)}2$ copies of the blow up of $\bP^1\times\bP^1\times\bP^1$ at some point.
\end{lemma}
\begin{figure}[h]
\centering
\includegraphics[height=2.3cm]{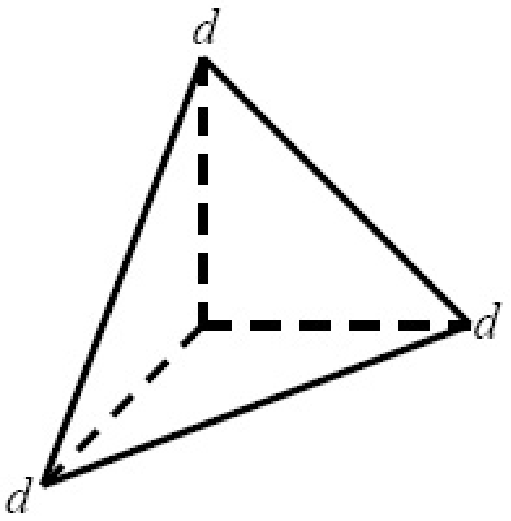}\caption{\small{$\D_d$ }}
\label{tetraedro1}
\end{figure}
\FloatBarrier
\noindent\textit{Proof}. Let $\D_d$ be the polytope defining $V_{3,d}$ as a toric variety; it
is a tetrahedron of side length $d$, as in Figure \ref{tetraedro1}. We start cutting $\D_d$ horizontally, so we get a $\D_{d-1}$ leaning
on one layer of height $1$, as in Figure \ref{tetrasezione}. We
denote by $S^1_k$ a layer of height $1$ and base of side length
$k$.
Iterating the cut in $\D_{d-1}$ we get a subdivision of $\D_d$ in
$(d-1)$ layers of height $1$ plus a tetrahedron on the top. Then we further subdivide each layer into a certain number of
polytopes, i.e. inside the layer $S_k^1$ there will be
$\frac{(k-1)(k-2)}{2}$ cubes, $k$ tetrahedra $\D_1$, and $k-1$
blocks $\Sigma$ as the one represented in Figure \ref{blocco}. The
polytope $\Sigma$ can be interpreted as obtained by removing a
tetrahedron from a cube; this corresponds to the blow up of $\bP^1\times\bP^1\times\bP^1$ at one point, which we denote
$\widetilde{(\bP^1)^3}$.
So far we obtained a subdivision $\G$ of $\D_d$ into cubes,
tetrahedra $\D_1$ and blocks $\Sigma$. Let $C_{ijk}$ be the cube having vertices
$(i,j,k),(i+1,j,k),(i,j+1,k),(i,j,k+1),(i+1,j+1,k),(i+1,j,k+1),(i,j+1,k+1),(i+1,j+1,k+1),$
with $i,j,k\in \mathbb{N}$. Our
piecewise linear function $F$ is defined as follows: on the cube
$C_{ijk}$ in the partition of $\D_d$, the function is given by the linear form
$$F\mid_{
C_{ijk}}(x_1,x_2,x_3)=(1+2i)x_1+(1+2j)x_2+(1+2k)x_3-(i+j+k+i^2+j^2+k^2).$$
An easy calculation shows that the
function is strictly convex on our subdivision $\G$. Hence the subdivision $\G$ is regular, and it yields a degeneration of $V_{3,d}$ to $\left({d \atop 3}\right)$ copies of $\bP^1\times\bP^1\times\bP^1$, corresponding to the cubes, plus $\frac{d(d+1)}2$ copies of
$\bP^3$, corresponding to the tetrahedra, and $\frac{d(d-1)}2$ copies of the blow up of $\bP^1\times\bP^1\times\bP^1$ at some point, corresponding to the blocks $\Si$.

\rightline{$\Box$}

\noindent In the following we are going to prove Alexander--Hirschowitz's Theorem
treating in separated sections the cases $d$ odd and $d$ even. We start with the following Lemmata, which will be useful throughout.

\begin{figure}[h]
\centering
\includegraphics[height=1.5cm]{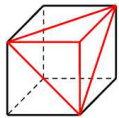}
\caption{\small{One of the $8$ tetrahedra $\D_1$ containing $3$ edges of $C$ meeting at one vertex.}}\label{tetra_tg}
\end{figure}

\begin{lemma}\label{tetraedro_tg}
Let $C$ be the unitary cube in $\bR^3$, and let $T$ be one of the $8$ tetrahedra $\D_1$ containing $3$ edges of $C$ meeting at one vertex. Then the toric variety it defines is a tangent $\bP^3$ to $(\bP^1)^3$ at one coordinate point of $\bP^7$.
\end{lemma}
\noindent\textit{Proof}. Note that $T$ is a tetrahedron like the one in Figure \ref{tetra_tg}. We can assume that $C$ is the cube with the origin as a vertex, and that $T$ is the tetrahedron of vertices $\{(0,0,0),(1,0,0),(0,1,0),(0,0,1)\}$. Then, an immediate computation shows that $T$ is the tangent space to the Segre variety $(\bP^1)^3$ defined by $C$, at $[1:0:0:0:0:0:0:0]$, one of the coordinate points of the $\bP^7$ spanned by $(\bP^1)^3$.

\rightline{$\Box$}

\begin{lemma}\label{observations}The following facts hold:
\begin{itemize}
\item[(i)]$\bP^3$ is not $0-$defective.
\item[(ii)]The Segre threefold $\bP^1\times\bP^1\times\bP^1$ is not $1-$defective.
\end{itemize}
\end{lemma}
\noindent\textit{Proof}. Assertion (i) is obvious by Definition \ref{defect}, since $Sec_0(\bP^3)=\bP^3$. In order to prove (ii), we recall that the Segre embedding of $(\bP^1)^3$ in $\bP^7$ is represented by a cube as a toric variety.
\begin{figure}[h] \centering
\includegraphics[height=1.5cm]{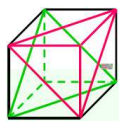}
\caption{\small{The subdivision of the cube.}}\label{disgiunti}
\end{figure}
If we look at Figure \ref{disgiunti}, we see that the two tetrahedra pointed out are of the type introduced in Lemma \ref{tetraedro_tg}, i.e. they correspond to tangent $\bP^3$'s to $(\bP^1)^3$. Since they are disjoint, the relative tangent $\bP^3$'s are skew in the ambient space $\bP^7$, so they span a linear space of maximal dimension $4\cdot 2-1=7$. If instead of the tangent $\bP^3$'s defined by the tetrahedra in the picture we take general ones, by semicontinuity we get that $\dim(Sec_1(\bP^1)^3)$ is maximal, hence $(\bP^1)^3$ is not $1-$defective.

\rightline{$\Box$}

\section{Case $d$ odd}

\noindent Let $d$ be odd, and let $\de$ be the degeneration of $V_{3,d}$ as in Lemma \ref{degeneration1}; by the same Lemma we have that each layer $S_k^1$ degenerates to $\frac{(k-1)(k-2)}{2}$ copies of $\bP^1\times\bP^1\times\bP^1$, $k$ copies of $\bP^3$ and $k-1$ copies of $\widetilde{(\bP^1)^3}$. We recall that a nonsingular polytope $\D$ describes a closed embedding of the toric variety $X(\D)$, i.e. it determines a very ample line bundle on $X(\D)$ \cite{bib:Ful}. Then $S_k^1$ corresponds to $\bP^3$ blown up at a point $p$, embedded in $\bP^{M_k-1}$, where $M_k=\frac{(k+2)(k+1)}2+\frac{(k+1)k}2$, via the proper transform of the linear system of surfaces of degree $k$ with a point of multiplicity $k-1$ at $p$. If we denote by $\pi$ the blow up,
the line bundle which embeds $\widetilde{\bP^3}$ in $\bP^{M_k-1}$ is
$\pi^*(\mathcal{O}_{\bP^3}(k))\otimes\mathcal{O}_{\widetilde{\bP^3}}(-(k-1)E)$
, where $E$ is the exceptional divisor of the blow up. We have the following:

\begin{lemma}\label{strati}
Let $k$ be odd; then $\widetilde{\bP^3}$, embedded in $\bP^{M_k-1}$, is not $\left(\frac{(k+1)^2}4-1\right)$--defective.
\end{lemma}

\begin{wrapfigure}[11]{r}[0pt]{4cm}
\includegraphics[height=3.5cm]{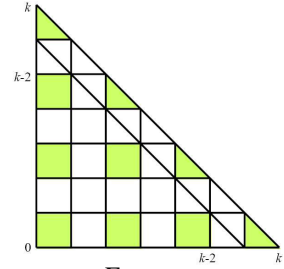}
\caption{\small{Projection of $S_k^1$}}\label{odd}
\end{wrapfigure}
\noindent\textit{Proof}. Consider the subdivision $\G_k$ of $S_k^1$ in cubes, $\Sigma$'s and $\D_1$'s. We have seen that it is regular; the unitary cube corresponds to the toric variety $(\bP^1)^3$ which is not $1$-defective, whereas the tetrahedron $\D_1$ defines $\bP^3$ as a not $0$--defective variety (see Lemma \ref{observations}). If we take disjoint polytopes as in Figure \ref{odd}, we can apply Proposition \ref{prop}. Indeed, in the layer $S_k^1$ we can point out $\frac{k^2-1}8$ disjoint cubes, and $\frac{k+1}2$ disjoint tetrahedra $\D_1$. Thus applying Proposition \ref{prop} we get that $\widetilde{\bP^3}$ in its embedding in $\bP^{M_k-1}$ is not $\left(2\dfrac{k^2-1}8+\dfrac{k+1}2-1\right)$--defective.
\bigskip

\noindent We are now ready to state:

\begin{lemma}\label{d odd}
Let $d$ be odd, then $V_{3,d}$ is not $n_d$--defective.
\end{lemma}
\noindent\textit{Proof}. Consider the subdivision of $\D_d$ into $d-1$ layers $S_k^1$ plus one $\D_1$. This subdivision is regular by an argument analogous to the one used in Lemma \ref{degeneration1}. By the previous Lemma we know that $S_k^1$ for $k$ odd is not $\left(\frac{(k+1)^2}4-1\right)$--defective. Let us take all the layers $S_k^1$ for $k\geq 3$ odd, and the tetrahedron $\D_1$ on the top: they are disjoint, therefore applying Proposition \ref{prop}, by induction on $d$ we get that $V_{3,d}$ is not $n_d$--defective. Indeed, let $d=5$: it is easy to see that $V_{3,5}$ is not $13$--defective, since $S_5^1$ is not $8$--defective, $S_3^1$ is not $3$--defective and $\D_1$ is not $0$--defective. Suppose now that we have proved that $V_{3,d-2}$ is not $n_{d-2}$--defective. We obtain the tetrahedron $\D_d$ adding to $\D_{d-2}$ the layers $S^1_{d-1}$ and $S^1_d$. So we just have to verify that $n_d+1=n_{d-2}+1+\frac{(d+1)^2}4$, which is an easy computation.

\rightline{$\Box$}

\begin{remark} {\rm To see that $V_{3,d}$ for $d$ odd is not $k$--defective
for any $k\in\mathbb{N}$, we just have to apply Lemma \ref{reduce} and Remark \ref{reduceremark},
since we have shown that $V_{3,d}$ is not $n_d$--defective.}
\end{remark}

\section{Case $d$ even}
\noindent When $d$ is even we cannot apply the same argument as for $d$ odd because the situation is more complicated; by Lemma \ref{expdim} we have to treat the cases $d=6,8,10,12$. Henceforth we will often regard the polytopes and the toric varieties they determine as the same objects.

\begin{figure}
\centering
\includegraphics[height=14cm]{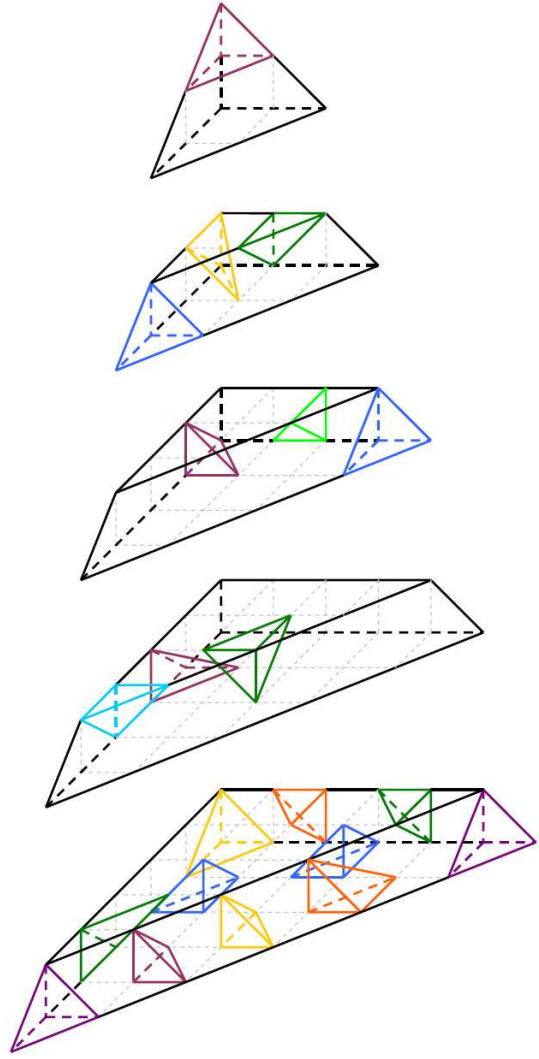}
\caption{\small{$\D_6$: it is decomposed in layers of height $1$ in which we have pointed out pairwise disjoint tetrahedra.}}\label{D6}
\end{figure}

\begin{lemma}\label{V6}
$V_{3,6}$ is not $20$--defective.
\end{lemma}
\noindent\textit{Proof}. Looking at Figure \ref{D6} one sees that among
the $21$ tetrahedra that we colored, $5$ are not of the type introduced
in Lemma \ref{tetraedro_tg}; we denote them by $\Si_1,\ldots,\Si_5$. They correspond to $\bP^3$ as toric varieties, but they are not
necessarily tangent $\bP^3$'s. It is clear that the $21$ $\bP^3$'s are independent, since their polytopes are disjoint, so they span a projective space of maximal dimension, $83$, that we denote by $\bP^M$. The four vertices of each tetrahedron correspond to independent points in the ambient space, so we can think of $\bP^M$ as the span of the $4$--tuples of vertices of the $21$ tetrahedra.
Let us observe that a
projectivity of $\bP^M$ can be represented by a $(84\times 84)$ matrix
up to multiplication by scalars. Fix $i=1,\ldots,5$, for example $i=1$,
and consider $\Si_1$; we have to distinguish two cases.
\begin{itemize}
 \item[Case I:] the tetrahedron $\Si_1$ lies inside a cube.
 We define $\f$ to be the
projectivity associated to a diagonal matrix $M_\f$, whose diagonal
is the string $(1,1,1,1,t,\ldots,t)$, where the first $4$ columns of
$M_\f$ correspond to the $4$ points generating $\Si_1$, and the remaining columns, $4$ by $4$ correspond to the ordered $4$-tuples of points spanning the other $\bP^3$'s.
 As we know $\Si_1$ lies inside a cube, which corresponds to a $(\bP^1)^3$; let $H_1$ be a general $\bP^3$ tangent to the $(\bP^1)^3$ where $\Si_1$ sits, at some general point. Since $(\bP^1)^3$ spans a $\bP^7$, we can choose points $p_1,\ldots,p_4\in H_1$ such that
$\langle p_1,\ldots,p_4\rangle=H_1$, and such that $p_i=[a_{i0}:\ldots:a_{i7}:0:\ldots:0]$ for any $i=1,\ldots,4$. Applying now $\f$ to $H_1$,
we get:
\begin{equation}\label{matrici}
\f\left(
  \begin{array}{cccc}
    a_{10} & a_{20} & a_{30} & a_{40}\\
    {\cdot \atop  \cdot} & {\cdot \atop \cdot}
     & {\cdot \atop \cdot} & {\cdot \atop \cdot} \\
    a_{17} & a_{27} & a_{37} & a_{47}\\
    0 & 0 & 0 & 0 \\
    {\cdot \atop \cdot} & {\cdot \atop \cdot}
     & {\cdot \atop \cdot} & {\cdot \atop \cdot} \\
    0 & 0 & 0 & 0
  \end{array}
\right)=
\left(
  \begin{array}{cccc}
    a_{10} & a_{20} & a_{30} & a_{40}\\
    {\cdot \atop  \cdot} & {\cdot \atop \cdot}
     & {\cdot \atop \cdot} & {\cdot \atop \cdot} \\
    a_{13} & a_{23} & a_{33} & a_{43}\\
    ta_{14} & ta_{24} & ta_{34} & ta_{44}\\
    {\cdot \atop  \cdot} & {\cdot \atop \cdot}
     & {\cdot \atop \cdot} & {\cdot \atop \cdot} \\
    ta_{17} & ta_{27} & ta_{37} & ta_{47}\\
    0 & 0 & 0 & 0 \\
    {\cdot \atop  \cdot} & {\cdot \atop \cdot}
     & {\cdot \atop \cdot} & {\cdot \atop \cdot} \\
    0 & 0 & 0 & 0
  \end{array}
\right)
\end{equation}
Taking the limit as $t\rightarrow 0$ we obtain that the limit image
of $H_1$ via $\f$ is the space $\mathscr{H}_1$ defined by the following matrix
$$\left(
  \begin{array}{cccc}
    a_{10} & a_{20} & a_{30} & a_{40}\\
     {\cdot \atop  \cdot} & {\cdot \atop \cdot}
     & {\cdot \atop \cdot} & {\cdot \atop \cdot} \\
    a_{13} & a_{23} & a_{33} & a_{43}\\
    0 & 0 & 0 & 0 \\
     {\cdot \atop  \cdot} & {\cdot \atop \cdot}
     & {\cdot \atop \cdot} & {\cdot \atop \cdot} \\
    0 & 0 & 0 & 0
  \end{array}
\right)$$
We note that
$\mathscr{H}_1$ is the image of $H_1$ via the
projection $\pi:\bP^7\dashrightarrow\bP^3$, such that
$\pi(x_0:\ldots :x_7)\mapsto (x_0:\ldots:x_3)$, which is the projection from the $\bP^3$ spanning with $\Si_1$ the $\bP^7$ where $(\bP^1)^3$ sits, i.e., the $\bP^3$ generated by the $4$ vertices of the cube not belonging to $\Si_1$. An easy calculation shows that the general $\bP^3$ tangent to $(\bP^1)^3$, in particular $H_1$, doesn't intersect the center of the projection.
Hence $\mathscr{H}_1$ is still a $\bP^3$, and it is precisely $\Si_1$.
\item [Case II:]The tetrahedron $\Si_1$ lies inside a block $\Si$ as in Figure \ref{blocco}, so the $\bP^3$ it defines sits in a $\bP^6$, since the toric variety defined by $\Si$ is a $(\bP^1)^3$ blown up
at a point. Let $\psi$ be
the projectivity of $\bP^M$ defined by a diagonal matrix $M_\psi$,
whose diagonal is the string $(1,1,1,1,t,\ldots,t)$, as well as in
Case I, where the first four columns correspond to the four points generating $\Si_1$. Let $H_1$ be a general $\bP^3$ tangent to the $\widetilde{(\bP^1)^3}$ where the $\bP^3$ corresponding to $\Si_1$ sits; choose $p_1,\ldots,p_4\in H_1$ such that
$\langle p_1,\ldots,p_4\rangle=H_1$. The space $H_1$ is tangent to $\widetilde{(\bP^1)^3}$, which spans a $\bP^6$, hence we can take the points $p_1,\ldots,p_4$
such that $p_i=[a_{i0}:\ldots:a_{i6}:0:\ldots:0]$ for any $i=1,\ldots,4$. In this setting we use the argument of Case I, and obtain that the limit of $H_1$ via $\psi$ is its image via the
projection $\pi:\bP^6\dashrightarrow\bP^3$, such that
$\pi(x_0:\ldots :x_6)\mapsto (x_0:\ldots:x_3)$, which is the projection
from the $\bP^2$ spanning with $\Si_1$ the $\bP^6$ where the blow up of
$(\bP^1)^3$ sits, i.e. the $\bP^2$ generated by the $3$ vertices of the block $\Si$ not belonging to $\Si_1$. An easy calculation shows that the general tangent $\bP^3$ doesn't intersect the center of projection, hence the limit of $H_1$ is still a $\bP^3$, and it is exactly the one determined by $\Si_1$.
\end{itemize}

\noindent Let us denote by $\Si_5$ the tetrahedron in Figure \ref{D6} that we see in the base layer $S^1_6$, lying in a block $\Si$. The other tetrahedra, $\Si_1,\ldots,\Si_4$, lie inside unitary cubes. We apply the construction of Case I to $\Si_1,\ldots,\Si_4$, and the construction of Case II to $\Si_5$. Composing the degenerations one by one, we obtain
Figure \ref{D6} as our limit configuration. Now, since the limiting $\bP^3$'s span a projective space of dimension $(4\cdot 21-1)$, hence, by semicontinuity of the dimension
we get a lower bound on $dim(Sec_{20}(V_{3,6}))$, and we can conclude that
$V_{3,6}$ is not $20$--defective.

\rightline{$\Box$}

\begin{lemma}\label{base8}
$(\widetilde{\bP^3})$ embedded in $\bP^{80}$ via
$\pi^*(\mathcal{O}_{\bP^3}(8))\otimes\mathcal{O}_{\widetilde{\bP^3}}(-7E)$
 is not $19$--defective.
\end{lemma}

\begin{figure}[h]
\centering
\includegraphics[height=4.5cm]{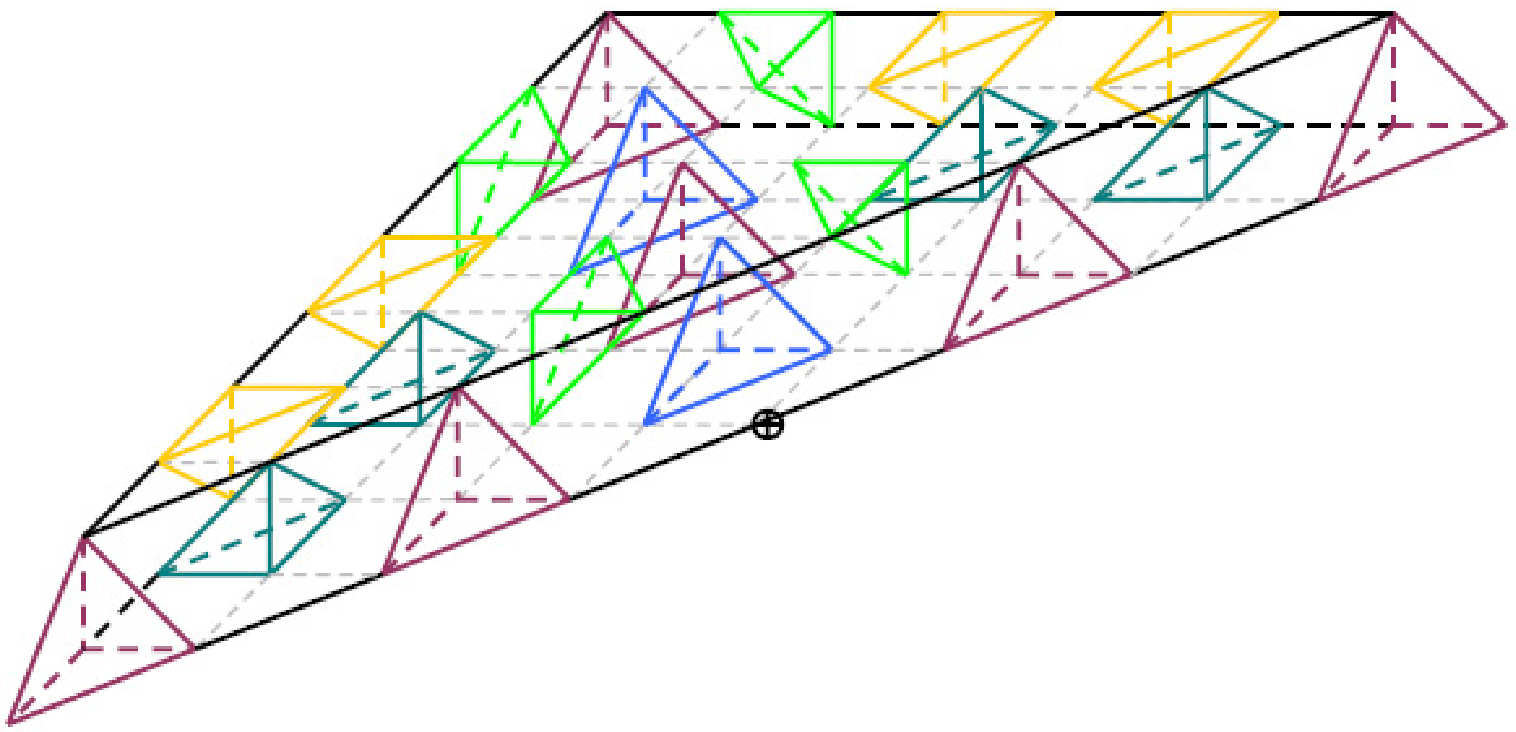}
\caption{\small{$S^1_8$}}\label{D8}
\end{figure}

\noindent\textit{Proof.} Figure \ref{D8} represents a layer $S^1_8$,
whose corresponding toric variety is $(\widetilde{\bP^3})$ in its
embedding in $\bP^{80}$. In the polytope we have colored $20$
tetrahedra of the type introduced in Lemma \ref{tetraedro_tg}.
Then the toric varieties they define are tangent $\bP^3$'s.
Since the tetrahedra are pairwise disjoint, the relative tangent
$\bP^3$'s are skew in the ambient space $\bP^{80}$, so they
span a linear space of maximal dimension $4\cdot 20-1$.
Notice that in the picture there is one point that is
not included in any tetrahedron. If instead of the tangent
$\bP^3$'s defined by the tetrahedra in the picture we take
general ones, by semicontinuity we still have
that $dim(Sec_{19}(\widetilde{\bP^3}))$ is maximal.

\rightline{$\Box$}

\begin{lemma}\label{base10}
$(\widetilde{\bP^3})$ embedded in $\bP^{120}$ via
$\pi^*(\mathcal{O}_{\bP^3}(10))\otimes\mathcal{O}_{\widetilde{\bP^3}}(-9E)$
 is not $29$--defective.
\end{lemma}

\begin{figure}[h]
\centering
\includegraphics[height=5.5cm]{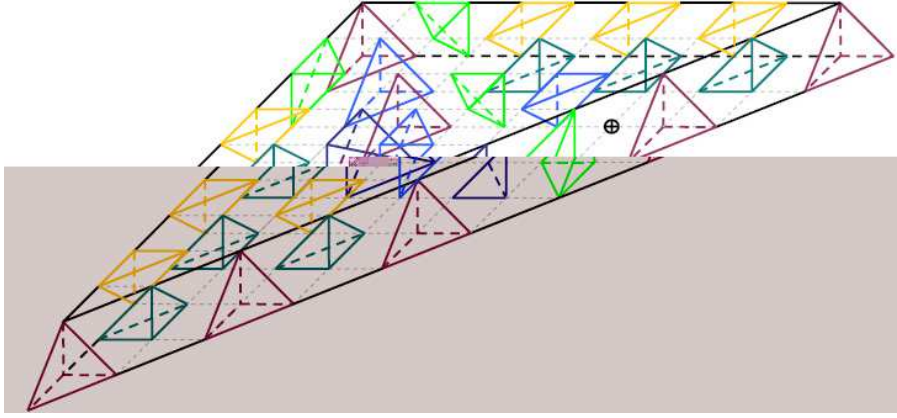}
\caption{\small{The layer $S^1_{10}$: the two tetrahedra in blue do not correspond to tangent $\bP^3$'s. }}\label{D10}
\end{figure}

\noindent\textit{Proof.} In Figure \ref{D10} we see a layer $S^1_{10}$
 which corresponds to the embedding of $(\widetilde{\bP^3})$ in $\bP^{120}$.
  The colored tetrahedra inside $S^1_{10}$ are $30$, and only two of them
  are not of the type described in Lemma \ref{tetraedro_tg}. Applying the
  argument of Lemma \ref{V6}, we regard Figure \ref{D10} as the limit of
  a configuration composed of $30$ independent tangent $\bP^3$'s at
  general points of $(\widetilde{\bP^3})$. The $\bP^3$'s of the limiting
  configuration span a projective space of dimension $(4\cdot 30-1)$.
  Then, again by semicontinuity of the dimension, we have that $Sec_{29}(\widetilde{\bP^3})$
   is not $29$--defective. As well as in $S^1_8$ we have one point left
   out of the configuration.

\rightline{$\Box$}
\begin{cor}\label{cor8-10}
$V_{3,8}$ and $V_{3,10}$ are not $k$--defective for any $k\in\mathbb{N}$.
\end{cor}

\noindent\textit{Proof.} Let us consider the configuration composed of
the tetrahedron of Figure \ref{D6} leaning on a layer $S^1_7$, which
in turn lies on a layer $S^1_8$ as the one in Figure \ref{D8}.
We denote by $D_8$ such a configuration; it can be seen as the
central fiber (composed of $\D_6$, $S^1_7$ and $S^1_8$) of a degeneration
of $V_{3,8}$. By Lemma \ref{V6} and Lemma \ref{base8}, we can apply
Proposition \ref{prop}, being $\D_6$ and $S^1_8$ disjoint, thus we
obtain that $dim(Sec_{n_8}(V_{3,8}))=4\cdot 41-1=163$, i.e. $V_{3,8}$
is not $40$--defective. Note that $Sec_{n_8}(V_{3,8})\subset\bP^{164}=\bP^{N_8}$,
hence $Sec_{n_8}(V_{3,8})$ is a hypersurface. We look now at the defectivity of $Sec_{n_8+1}(V_{3,8})$; let us denote by $H_{n_8}$ the general tangent space to
$Sec_{n_8}(V_{3,8})$, which is a hyperplane, and by $H_{n_8+1}$
the tangent space to $Sec_{n_8+1}(V_{3,8})$ at a general point.
Then, by Terracini's Lemma we get that $H_{n_8+1}=\langle H_{n_8},T_p(V_{3,8})\rangle$,
$p$ being a general point on $V_{3,8}$; hence $dim(H_{n_8+1})=N_8$.
It follows that $V_{3,8}$ is not $k$--defective for any $k>n_8$, hence
for any $k\in\mathbb{N}$ (see Lemma \ref{reduce}). We use a similar
argument to show that $V_{3,10}$ is not $n_{10}$--defective: let us consider
the configuration $D_{10}$ composed of $D_8$ leaning on two layers,
one of type $S^1_9$, below which we put a $S^1_{10}$ as the one in
Figure \ref{D10}. Since $D_8$ and $S^1_{10}$ are
disjoint polytopes, we apply again Proposition \ref{prop} to this configuration,
using Lemma \ref{base10}, thus we get that $V_{3,10}$ is not $70$--defective. Now we want
to extend the non $k$--defectivity of $V_{3,10}$ to any $k\in\mathbb{N}$.
If $k<n_{10}$ we apply Lemma \ref{reduce}; let us now consider $Sec_{n_{10}+1}(V_{3,10})$, we call
$H_{n_{10}}$ the general tangent space to $Sec_{n_{10}}(V_{3,10})$, and
$H_{n_{10}+1}$ the tangent space to $Sec_{n_{10}+1}(V_{3,10})$ at some general point.
Looking at the polytopes which compose $D_{10}$, we see that two vertices
stay out of the configuration (one pointed out in $D_8$, the other one in $S^1_{10}$);
indeed $dim(Sec_{n_{10}}(V_{3,10}))=N_{10}-2$. Hence $codim(H_{n_{10}})=2$.
Again by Terracini's lemma we get that $H_{n_{10}+1}=\langle H_{n_{10}},T_p(V_{3,10})\rangle$,
with $p$ a general point on $V_{3,10}$. Then $dim(H_{n_{10}+1})=N_{10}$;
if it were not, then $H_{n_{10}+1}$ would be a hyperplane in $\bP^{N_{10}}$;
let $\pi_{n_{10}}$ be the projection from $H_{n_{10}}$ in $\bP^{N_{10}}$ onto $\bP^1$.
Since $T_p(V_{3,10})\subset H_{n_{10+1}}$, its projection is a point.
But $\pi_{n_{10}}(T_p(V_{3,10}))$ is exactly the general tangent space to the
projection of $V_{3,10}$, which must be a point. Hence $V_{3,10}$ would be degenerate,
which is a contradiction. With this argument we get that $Sec_k(V_{3,10})$
is not $k$--defective for any $k>n_{10}$.

\rightline{$\Box$}
\begin{figure}[h]
\centering
\includegraphics[height=15cm]{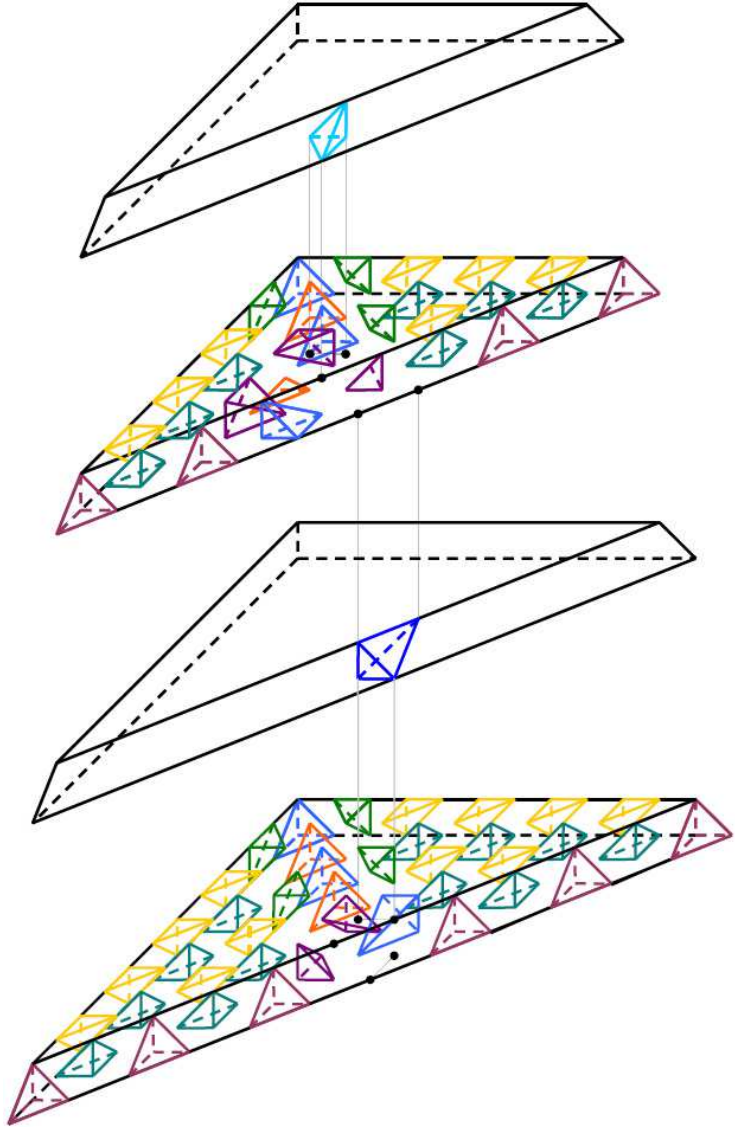}
\caption{\small{The layers $S^1_{9},S^1_{10},S^1_{11},S^1_{12}$ in $D_{12}$}}\label{D12}
\end{figure}


\begin{lemma}\label{lemma12}
$V_{3,12}$ is not $k$--defective for any $k\in\mathbb{N}$.
\end{lemma}

\noindent\textit{Proof.} Let us look at Figure \ref{D12};
it is a composition
of layers of type $S^1_9$, $S^1_{10}$, $S^1_{11}$ and $S^1_{12}$.
We place over the four layers represented, a configuration of type $D_8$,
as in Corollary \ref{cor8-10}. We call this new configuration $D_{12}$;
we are going to show that $V_{3,12}$ is not $n_{12}$--defective using the
strategy of Lemma \ref{V6}. Looking at the $113$ colored tetrahedra in
$D_{12}$ (counting also the ones in $D_8\subset D_{12}$),
most of them correspond to tangent $\bP^3$'s, as in Lemma \ref{tetraedro_tg},
but $13$ of them are not tangent $\bP^3$'s, so we want to think of them as
limits of general tangent $\bP^3$'s. If we denote by $\Si_1,\ldots,\Si_{13}$
the mentioned tetrahedra, we see that some of them lie inside unitary cubes,
but we find three tetrahedra, one lying in $D_6$ ($D_6\subset D_8$, moreover see Lemma \ref{V6}), one in $S^1_{10}$, the other one in
$S^1_{11}$, which sit inside blocks $\Si$ as in Figure \ref{blocco}. Let us call them $\Si_1$, $\Si_2$ and $\Si_3$. For all the tetrahedra $\Si_1,\ldots,\Si_{13}$ we can reapply the techniques used in Lemma \ref{V6}, i.e.
one by one we compose the degenerations of general tangent $\bP^3$'s to the limit $\bP^3$'s defined by
the relative $\Si_i$'s, $i=1,\ldots,13$. Of course we will apply the construction of Case II to $\Si_1,\Si_2$ and $\Si_3$, and the construction of Case I to $\Si_4,\ldots,\Si_{13}$.
At the end of the process, composing
all the degenerations, we obtain $D_{12}$ as our limit configuration.
Finally we observe that, since the $113$ limit $\bP^3$'s span a projective
space of dimension $(4\cdot 113-1)$, hence, by semicontinuity of the dimension
we get a lower bound on $dim(Sec_{n_{12}}(V_{3,12}))$, and we can conclude that
$V_{3,12}$ is not $n_{12}$--defective. But we are interested in the
$k$--defectivity of $V_{3,12}$ for any $k\in\mathbb{N}$. So, we are going
to study $Sec_{n_{12}+1}(V_{3,12})$. Let us notice that in the configuration
$D_{12}$ three points are not contained in any of the pointed out tetrahedra. However they are very close to each other, so we can think
of the three points as being contained in a further tetrahedron in $D_{12}$. Having in mind a configuration with $n_{12}+2$ disjoint tetrahedra, we can repeat the
argument we used before. Then
again by semicontinuity we obtain that $Sec_{n_{12}+1}(V_{3,12})$ spans
the whole ambient space. It follows that $V_{3,12}$ is not
$k$--defective for any $k>n_{12}$, hence for any $k$ (see Lemma \ref{reduce}).

\rightline{$\Box$}
\medskip

\noindent So far we have shown that $V_{3,d}$ is not $k$--defective for $d=6,8,10,12$
and any $k\in\mathbb{N}$. In the following we will study the cases $d\equiv6 \mod 8$,
$d\equiv0 \mod 8$, $d\equiv10 \mod 8$, which we will handle using recurrence on
the configurations we have already built.

\begin{figure}[h]
\centering
\includegraphics[height=4.5cm]{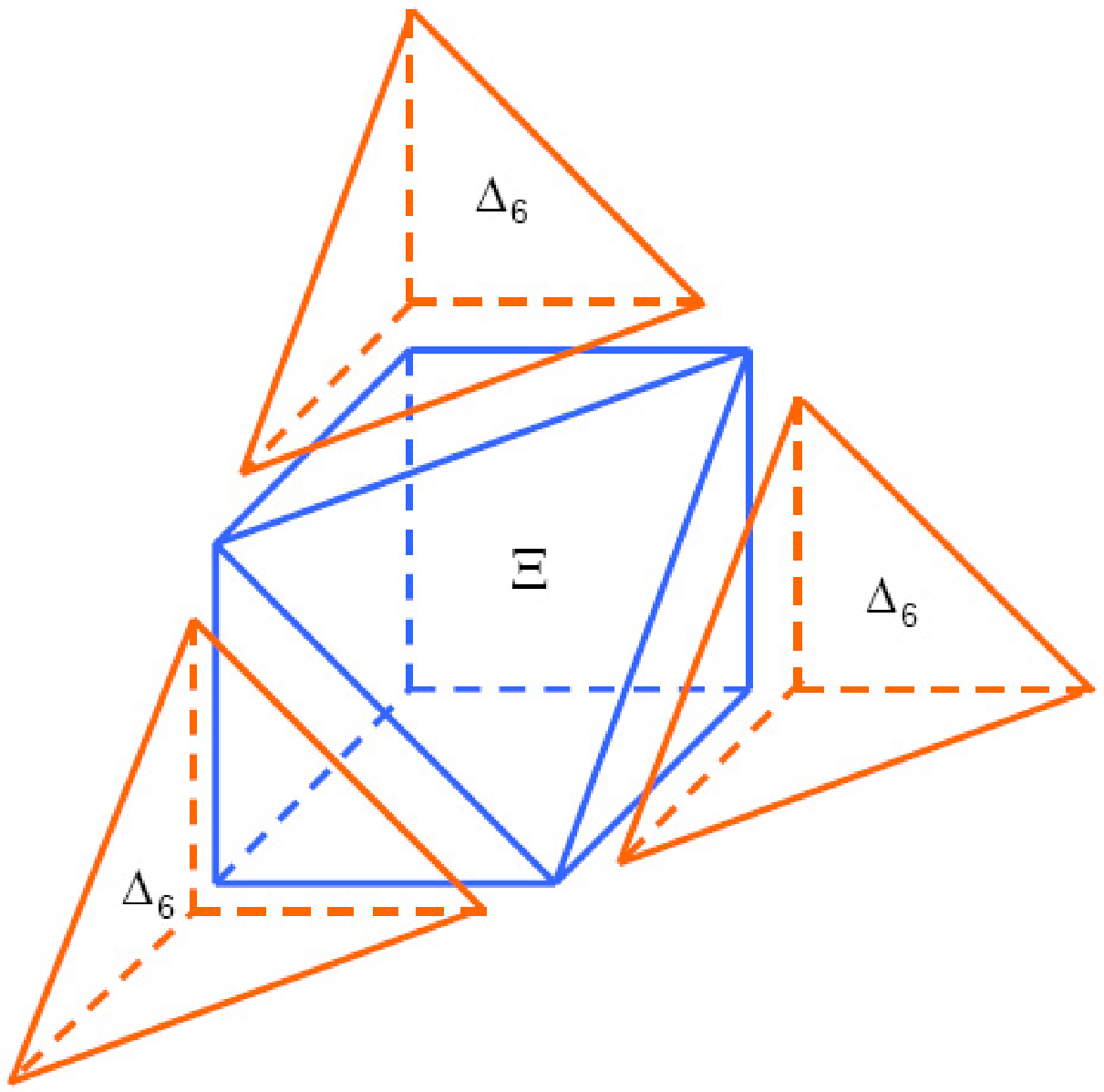}
\caption{\small{$\D_{14}$}}\label{d=14}
\end{figure}

\subsection{Case d\,=\,6+8k}\label{subsection1}

\begin{figure}[h]
\centering
\includegraphics[height=3cm]{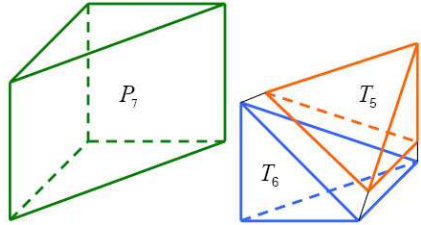}
\caption{\small{Subdivision of the block $\Xi$}}\label{base-d=14}
\end{figure}

\begin{figure}[h]
\centering \subfigure[Orthogonal projection of $P_7$.]
{\includegraphics[height=1.7cm]{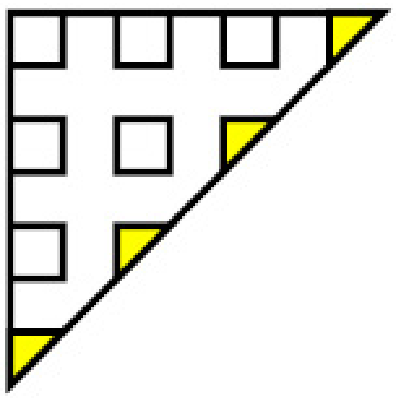}}\qquad\qquad\qquad
\qquad \subfigure[$\gamma_7$.]{\includegraphics[height=2.2cm]{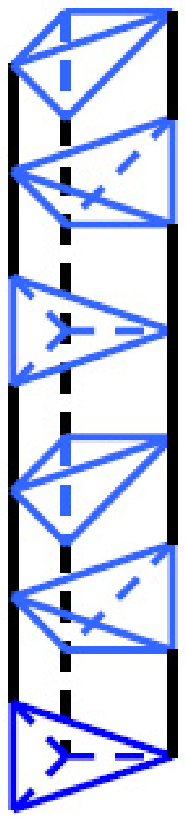}}
\caption{\small{}} \label{proiez.semicubo-d=14}
\end{figure}

\begin{figure}[h]
\centering
\includegraphics[height=4.5cm]{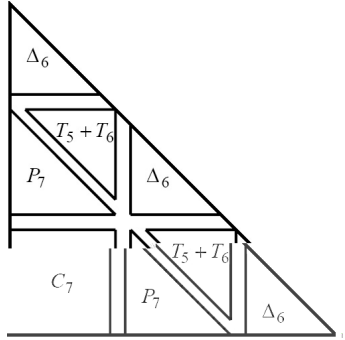}
\caption{\small{Orthogonal projection of $S_{22}^7$.}}\label{d=22}
\end{figure}

We start by examining the case $d=6+8k$. When $k=1$, the tetrahedron
$\D_{14}$  (see Lemma \ref{degeneration1}) can be
split up into one $\D_6$, that we triangulate following Figure \ref{D6},
and one layer $S_{14}^7$ of
height $7$ and big base of side length $14$. We can furthermore
divide $S_{14}^7$ in two triangulated tetrahedra $\D_6$ lying in its
ends, and a block $\Xi$, as illustrated in Figure \ref{d=14}. We can
triangulate $\Xi$ as follows: one semicube $P_7$ of side length $7$, plus one tetrahedron of side length $5$, that we denote by $T_5$,
and one tetrahedron of side length $6$, $T_6$. None of $T_5$ and
$T_6$ is a $3$-simplex, but we can think of them as the configurations
of $V_{3,5}$ and $V_{3,6}$ respectively (the subdivision of $\Xi$
is represented in Figure \ref{base-d=14}).
Recall that $V_{3,5}$ and $V_{3,6}$ are respectively not $13$--defective and
$20$--defective. In Figure \ref{proiez.semicubo-d=14}(a) we have represented the
orthogonal projection of the semicube $P_7$, which shows that it contains $6$
disjoint columns of height $7$, each containing $4$ pairwise disjoint
unitary cubes, and every cube is not $1$--defective; the squares in
Figure \ref{proiez.semicubo-d=14}(a)
are precisely the
projections of the columns. The yellow
triangles correspond to the projections of some different columns of
height
$7$, that we denote by $\gamma_7$, each containing $6$ pairwise
disjoint tetrahedra, as in Figure \ref{proiez.semicubo-d=14}(b). If we consider the toric
variety defined by $\gamma_7$, it is a $\bP^1\times\bP^2$ embedded in
$\bP^{23}$ via $p_2^*(\mathcal{O}_{\bP^2}(1))\otimes p_1^*(\mathcal{O}_{\bP^1}(7))$,
where $p_1,p_2$ are the projections on the factors; we see that it is not
$5$--defective. To show this, we will use the techniques of Lemma \ref{V6}.
We want to prove that the $6$ tetrahedra that we pointed out in $\gamma_7$
correspond to limits of general tangent $\bP^3$'s. Four of them are not
of the type introduced in Lemma \ref{tetraedro_tg}; let us denote them by
$\Si_1,\ldots,\Si_4$; every $\Si_i$ lies
inside a semicube, which corresponds to $\bP^1\times\bP^2$ and spans a $\bP^5$. So in order to use the argument of Lemma \ref{V6} we just have to pay attention to the ambient space. An easy calculation such as in Lemma \ref{V6} yields that $\Si_1,\ldots,\Si_4$ are limits of general tangent $\bP^3$'s hence
Figure \ref{proiez.semicubo-d=14}(b) is our limit configuration. Notice that the tetrahedra
corresponding
to the limit $\bP^3$'s are pairwise disjoint, then they span a projective
space of dimension $(4\cdot 6-1)$. Hence, by semicontinuity we can conclude
that the toric variety defined by $\gamma_7$ is not $5$--defective. Going
back to $P_7$, since it is composed of $4$ columns $\gamma_7$ plus $7\cdot 4$
unitary cubes and all these polytopes are pairwise disjoint, applying
Proposition \ref{prop} we get that the toric variety associated to $P_7$
is not $71$--defective. Recalling that above $S^7_{14}$ lies a $\D_6$,
which represents a $V_{3,6}$, and that in our configuration (see Figure \ref{d=14})
we have disjoint polytopes, then we can apply Proposition \ref{prop} and
obtain that $V_{3,14}$ is not $n_{14}$--defective (where $n_{14}+1=170$).
In order to study the generic case $d=6+8k$, we look at the tetrahedron
$\D_{6+8k}$; we proceed by recurrence. Suppose we already know that
$V_{3,d-8}$ is not $n_{d-8}$--defective. We build a configuration of
$\D_{d}$ in the following way: we set a layer of height $7$
with base of side length $d$, i.e. $S_d^7$, below $\D_{d-8}$, whose
configuration we already know by induction; now we subdivide $S_d^7$ as
follows: first we note that if $d=6+8k$,
then $S_d^7$ contains $\frac{k(k-1)}2$ cubes of side length $7$
($C_7$). In $C_7$ we insert $16$ disjoint columns of height $7$,
each containing $4$ disjoint unitary cubes, which, in turn, are not
$1$--defective, as shown in Lemma
\ref{observations}. Hence, being $16\times 4\times 2=128$,
we get that $C_7$ is not $127$--defective. Moreover,
$S_d^7$ contains $k$ semicubes equal to $P_7$, that we triangulated
before, and finally, on the skew stripe located in the front we set
$k+1$ copies of $\D_6$, $k$ copies of $T_5$ and $k$
ones of $T_6$, as in Figure \ref{d=22}, which represents the
orthogonal projection of $S_{22}^7$. We summarize our construction
in the following computation:
\begin{equation}\label{conti6+8k}
\D_d=\D_{d-8}+\left(\dfrac{k(k-1)}2\cdot C_7\right)+ (k\cdot P_7)+
(2\cdot\D_6)+ k\cdot (T_5+T_6)+ (k-1)\cdot\D_6
\end{equation}

\noindent Hence, recalling the defectivity of the toric varieties
associated to the polytopes composing $\D_d$, we observe that
$$n_{d}=n_{d-8}+\dfrac{k(k-1)}2\cdot 128+k\cdot 72+2\cdot 21+k\cdot (14+21)+(k-1)\cdot 21.$$
Now, noticing that the polytopes involved are
all disjoint, we can apply Proposition \ref{prop}, and we get that
$V_{3,6+8k}$ is not $n_{6+8k}$--defective.

\subsection{Configuration for d\,=\,8+8k}
\begin{figure}[h]
\includegraphics[height=1.5cm]{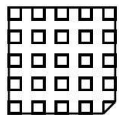}\caption{\small{Orthogonal projection of $H_9$.}}\label{H9}
\end{figure}

\begin{figure}[h]
\includegraphics[height=7.5cm]{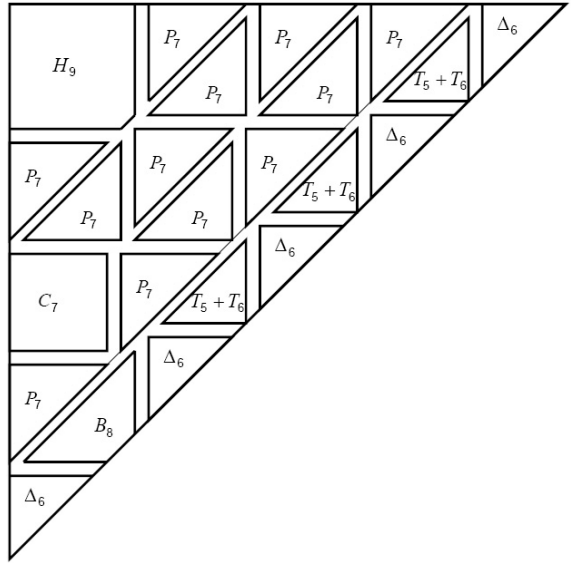}
\caption{\small{Orthogonal projection of $S^7_{40}$.}}\label{d=40}
\end{figure}

The case $d=8+8k$ with $k\geq 1$ is more complicated. As before, let
$k=1$. Then $d=16$, and we can decompose $\D_{16}$ setting a $\D_8$
above a $S_{16}^7$--layer.
\begin{figure}[h]
\centering
\includegraphics[height=4cm]{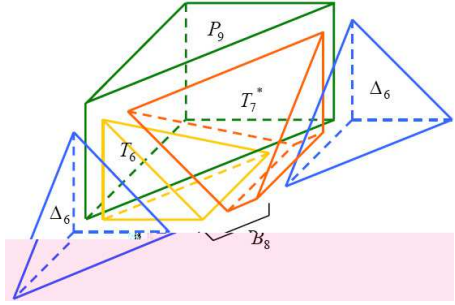}
\caption{\small{Subdivision of $S_{16}^7$.}}\label{d=16}
\end{figure}
\FloatBarrier
\noindent The problem is now to find a configuration for
$S_{16}^7$: we insert a prism having triangular basis of side length
$9$ and height $7$, denoted $P_9$;
it can be viewed as the half of a parallelepiped of height $7$
having square basis of side $9$, called $C_9$; then, on the skew
stripe in the front we put two tetrahedra $\D_6$ and a
block $B_8$ composed of a $T_6$ plus a cut tetrahedron of side $7$
which we denote $T_7^*$, as it is shown in Figure \ref{d=16}.
\begin{figure}[h] \centering
\subfigure[$T_7^*$ obtained from $\D_7$ removing the black prism.]
{\includegraphics[height=2.5cm]{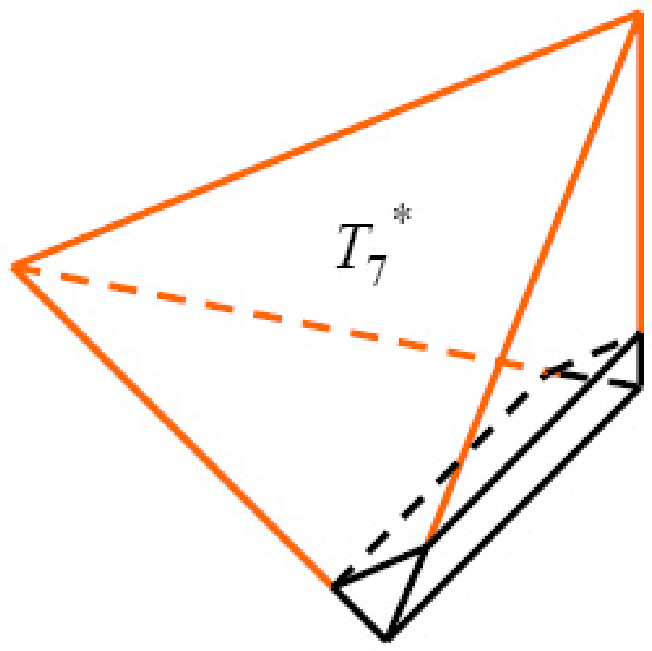}} \qquad\qquad\qquad
\subfigure[Tetrahedra in $T^*_7$.]
{\includegraphics[height=3.3cm]{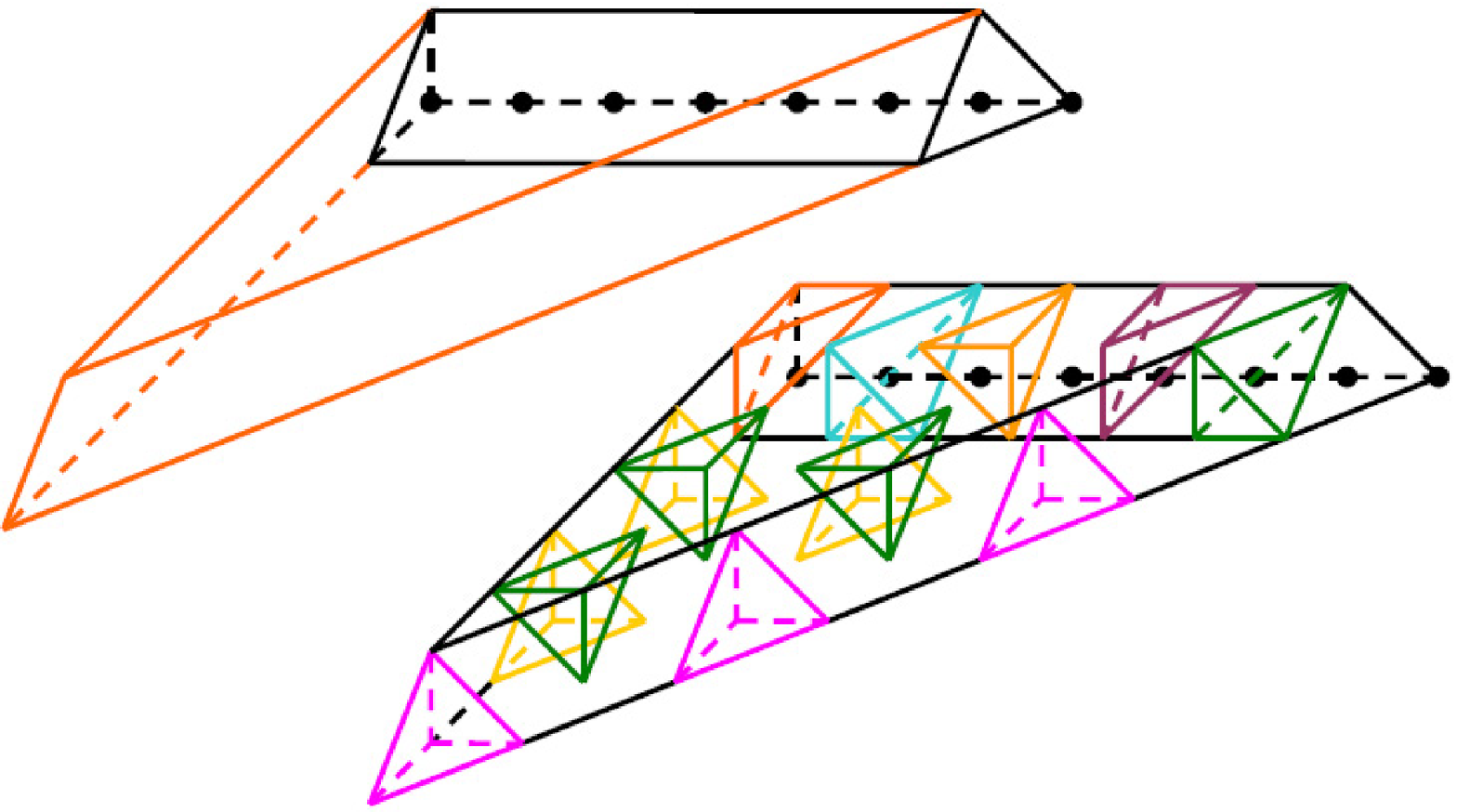}}\caption {\small{The
cut.}}\label{taglio7}
\end{figure}
\noindent We can obtain $T^*_7$ setting a $\D_5$ on a layer $S^1_6$,
in turn lying on a $S^1_7$, from which we remove the black prism as shown in
Figure \ref{taglio7}(a). We want to show that $T^*_7$ is not $27$--defective.
Notice that $\D_5$ is not $13$--defective; we focus on the cut $S^1_7$, which is
disjoint from $\D_5$. Let us consider the configuration in Figure
\ref{taglio7}(b). We see that $10$ tetrahedra correspond to the ones
introduced in Lemma \ref{tetraedro_tg}, i.e. they are tangent $\bP^3$'s.
But there are $4$ tetrahedra adjacent to the cut, which do not correspond
to tangent $\bP^3$'s. Each of them lies inside a semicube, so we apply the
argument of Lemmas \ref{V6} and \ref{lemma12}, that we also saw for the column $\gamma_7$, and we get that the tetrahedra
inside $S^1_7$ can be interpreted as limits of tangent $\bP^3$'s.
Being disjoint, they span a projective space of dimension $14\cdot 4-1$,
hence the cut layer $S^1_7$ is not $13$--defective (again with abuse of notation).
The block $B_8$ is not $48$--defective. Summing up, by Proposition \ref{prop}
we get that $V_{3,16}$ is not $n_{16}$--defective. In the
general case, when $d=8+8k$, we consider the parallelepiped $C_9$,
from which we remove a vertical edge of height $7$: in this way we
obtain a prism, $H_9$, that we can see in Figure \ref{H9}, and that
contains $96$ disjoint unitary cubes plus one column $\gamma_7$.
\noindent So $H_9$ is not $197$--defective. Just like
before we denote by $P_7$ the half of a cube $C_7$; $P_7$ is not
$71$--defective. We are now able to split $\D_d$ into $\D_{d-8}$,
and a $S^7_d$. In this layer we insert
(see Figure \ref{d=40}) in diagonal order: one $H_9$; $4$ copies of
$P_7$; then we have $\sum_{i=1}^{k-3}(i\cdot C_7+4\cdot
P_7)$; $k$ copies of $P_7$; in the front we put $(k+1)$ copies of
$\D_6$, $(k-1)$ copies of the block $T_5+T_6$, and finally one block
$B_8$.
\noindent Hence we have that:
$$\D_d=\D_{d-8}+H_9+4\cdot P_7+\sum_{i=1}^{k-3}(i\cdot C_7+4\cdot P_7)
+k\cdot P_7+(k-1)\cdot (T_5+T_6)+(k+1)\cdot \D_6+B_8.$$ Substituting
in the last expression the defectivity of each summand, we
get:$$n_{d}=n_{d-8}+198+4\cdot 72+\sum_{i=1}^{k-3}(i\cdot
128+4\cdot 72) +k\cdot 72+(k-1)\cdot (14+21)+(k+1)\cdot 21+49.$$
\bigskip
So applying Proposition \ref{prop} we get that $V_{3,8+8k}$ is not
$n_{8+8k}$--defective.

\subsection{Configuration for d\,=\,10+8k}
\begin{figure}[h] \centering
\subfigure[The cut of $S^3_9$ in $T_9^*$.]
{\includegraphics[height=3cm]{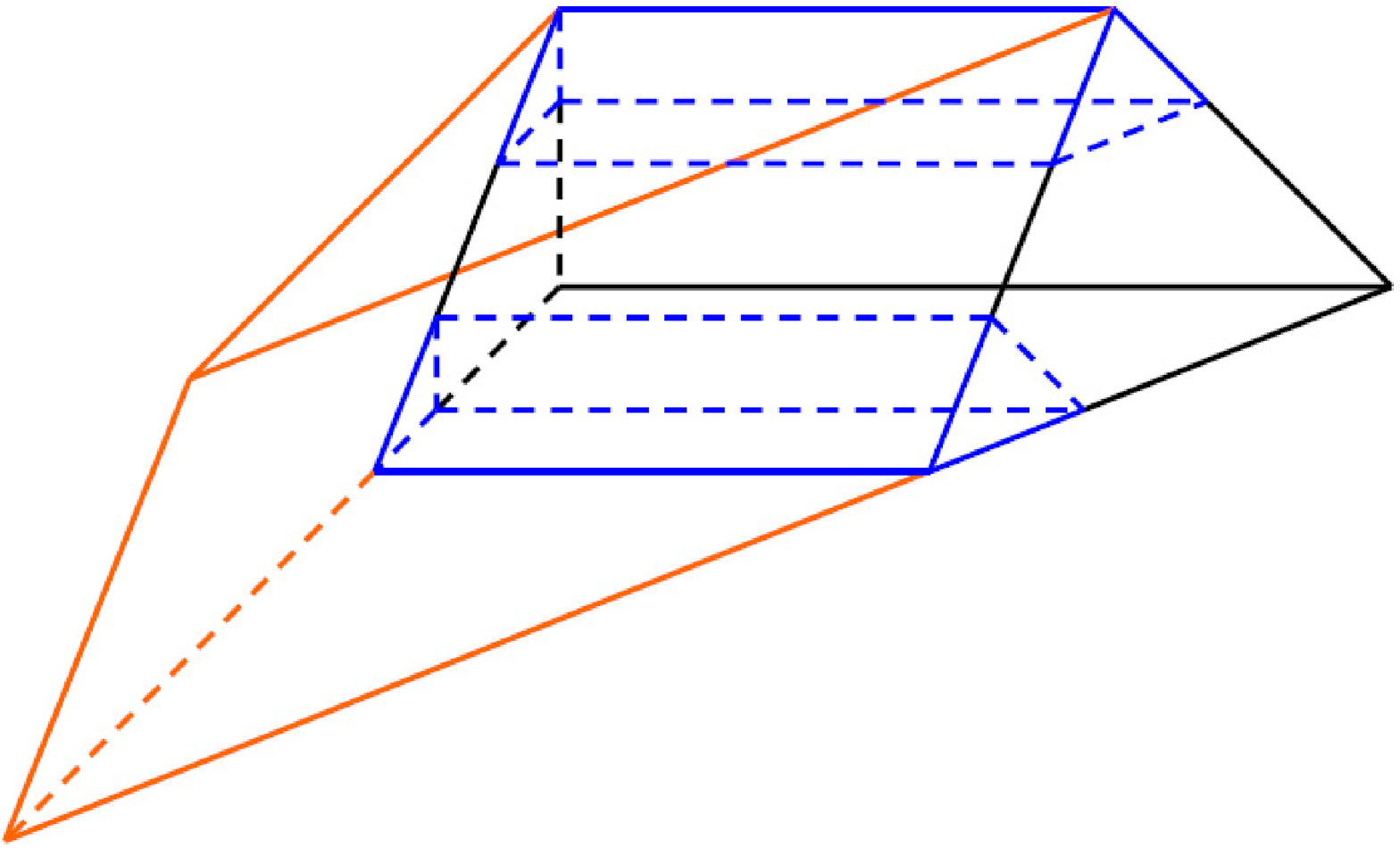}} \qquad\qquad
\subfigure[Orthogonal projection of $A_2$.]
{\includegraphics[height=2cm]{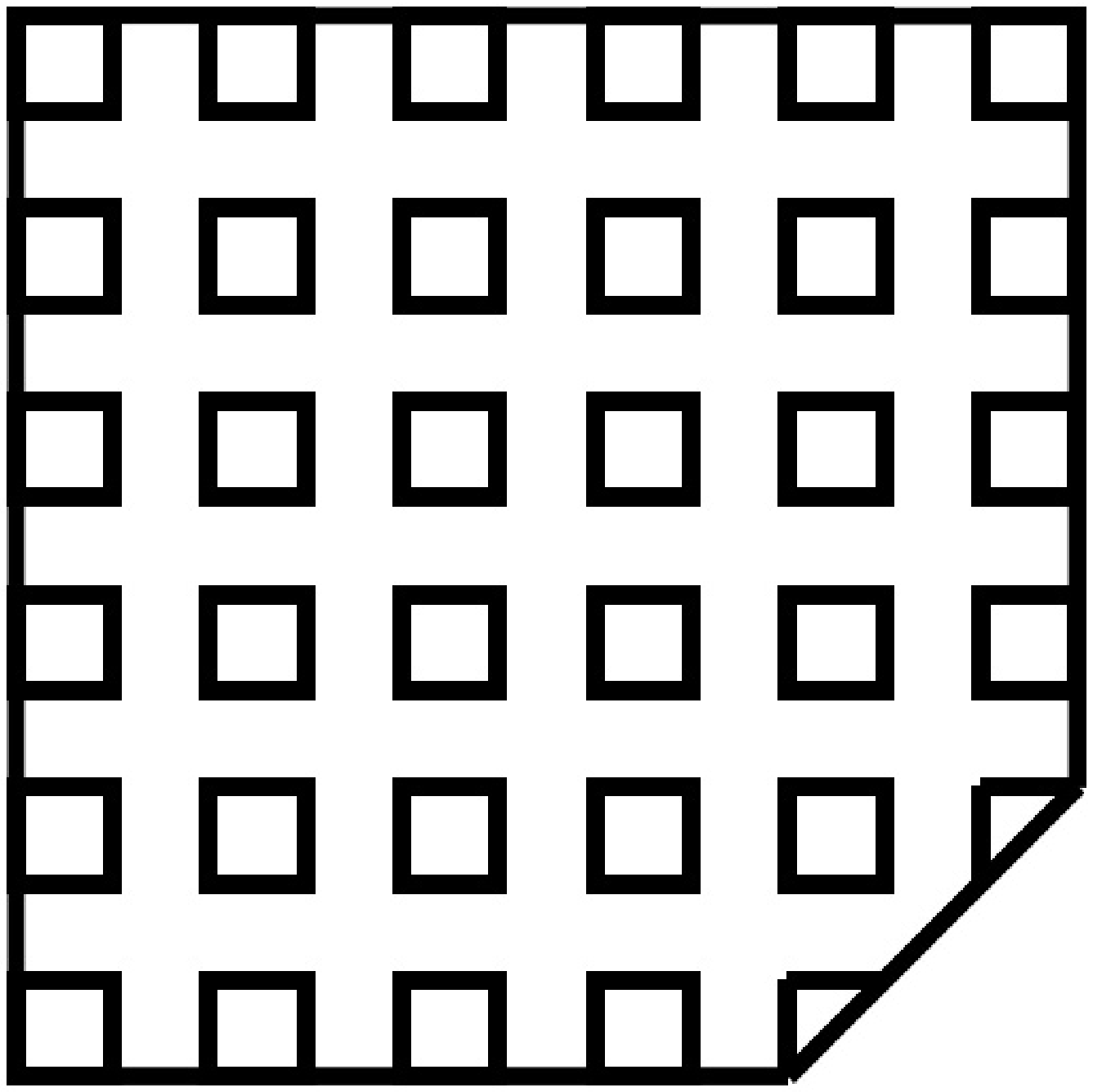}}\qquad\qquad
\subfigure[Orthogonal projection of $A_3$.]
{\includegraphics[height=3.2cm]{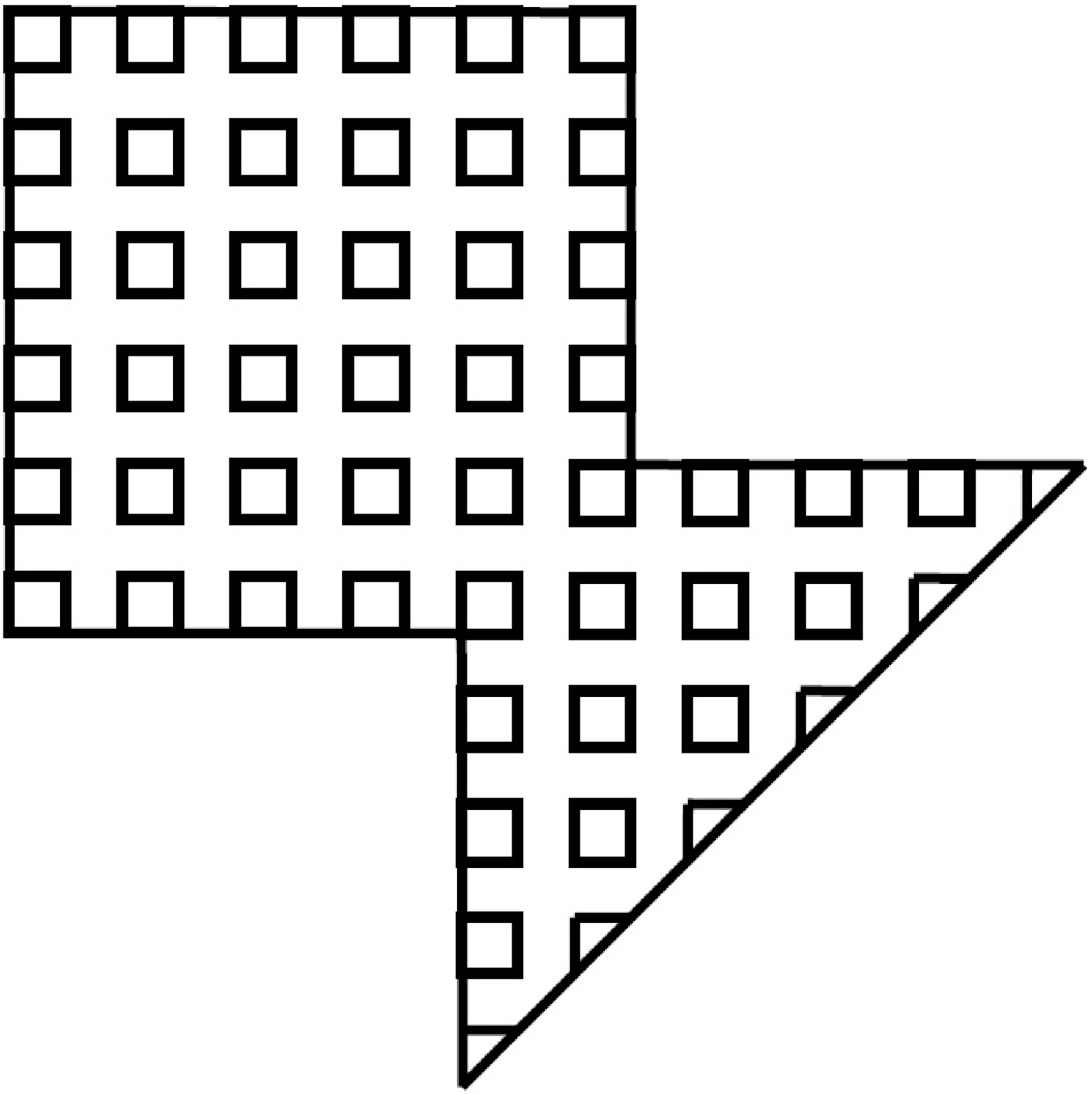}}\caption {\small{}}\label{taglio9}
\end{figure}

Let us examine the case $k=1$, corresponding to $d=18$; the
tetrahedron $\D_{18}$ is composed of a $\D_{10}$
leaning on a layer $S^7_{18}$.
\begin{figure}[h]
\centering
\includegraphics[height=4cm]{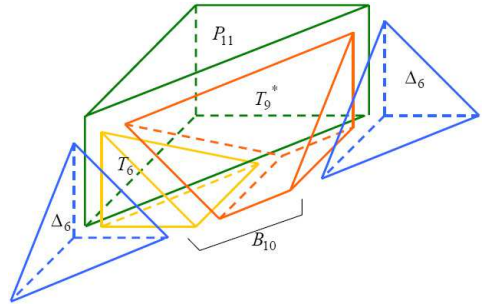}
\caption{\small{Subdivision of $S^7_{18}$.}}\label{d=18}
\end{figure}
\FloatBarrier \noindent Let us subdivide $S^7_{18}$ (see Figure
\ref{d=18}): we call $P_{11}$ the half of a parallelepiped of height
$7$ and basis of side length $11$; furthermore, in the front we set
two copies of $\D_6$ and a block $B_{10}$ composed of a $T_6$ and a
tetrahedron of side $9$ from which we remove a prism, as we did for
$T_7^*$. We obtain a solid, $T_9^*$; we want to show that it is not
$41$--defective. Let us decompose $T_9^*$ in the following way: a tetrahedron $\D_5$ leaning on a layer $S^1_6$, and this lies on a  $S^3_9$ that we cut as in Figure \ref{taglio9}(a). In particular we remove the part colored in black and blue, and in the remaining part we dispose the tetrahedra. Looking at Figure \ref{taglio9}(a) we see that $S^3_9$ with its cut is composed of $3$ layers, $S^1_7$, $S^1_8$ and $S^1_9$, all having a cut as well, so we dispose the tetrahedra in the $S^1_7$ and in $S^1_9$ as in Figure \ref{taglio7}(b); indeed the blue parts correspond to the black solid that we took off in Figures \ref{taglio7}(a) and \ref{taglio7}(b). Let us observe that the tetrahedra lying in $T_9^*$ are $42$.
\noindent Repeating the argument used for $T^*_7$, it
follows that the block $B_{10}$ is not $62$--defective. In $P_{11}$ we
insert $15$ disjoint columns of height
$7$ and square basis of side $1$, each containing $4$ disjoint
cubes; moreover $P_{11}$ contains $6$
columns $\gamma_7$. Then $P_{11}$ is not $155$--defective. Summing up
and using Proposition \ref{prop} we get that $V_{3,18}$ is not
$n_{18}$--defective. In the general
case $d=10+8k$, as usual we divide $\D_d$ in one triangulated
$\D_{d-8}$ plus one $S_{d}^7$. In the layer $S_{d}^7$ we insert: $2(k-2)$
copies of $C_7$, $2$ copies of $P_7$, one block $B_{10}$, $k+1$ copies
of $\D_6$, $k-1$ copies of $T_5+T_6$ and finally one block denoted
by $A_k$. Let us show the configuration of this block: for $k=2$, the block $A_2$ contains $33$ disjoint columns of height $7$ and square basis of
side $1$, and $2$ columns $\gamma_7$, hence it is not $275-$defective. The orthogonal projections
of these columns can be viewed in Figure \ref{taglio9}(b). For $k\geq 3$ the block $A_k$ is different: in Figure \ref{taglio9}(c) we see the projection of $A_3$. With this picture in mind, we observe that $A_k$ contains columns of height $7$ and side $1$ (not $7$--defective),
and columns $\gamma_7$ (not $5$--defective). In particular it contains
$\left(35+\sum_{i=3}^{4(k-2)+1}i\right)$ columns of the first type, and $(4\cdot (k-2)+2)$
columns $\gamma_7$. Thus if we set:
$\al_2+1:=276$, and for $k\geq 3$, $\al_k+1:=8\cdot\left(35+\sum_{i=3}^{4(k-2)+1}i\right)+6\cdot [4\cdot (k-2)+2]$, we have that
$A_k$ is not $\al_k$--defective for any $k\geq 2$.
\begin{figure}[h] \centering
\subfigure[Orthogonal projection of $S^7_{26}$.]
{\includegraphics[height=5.5cm]{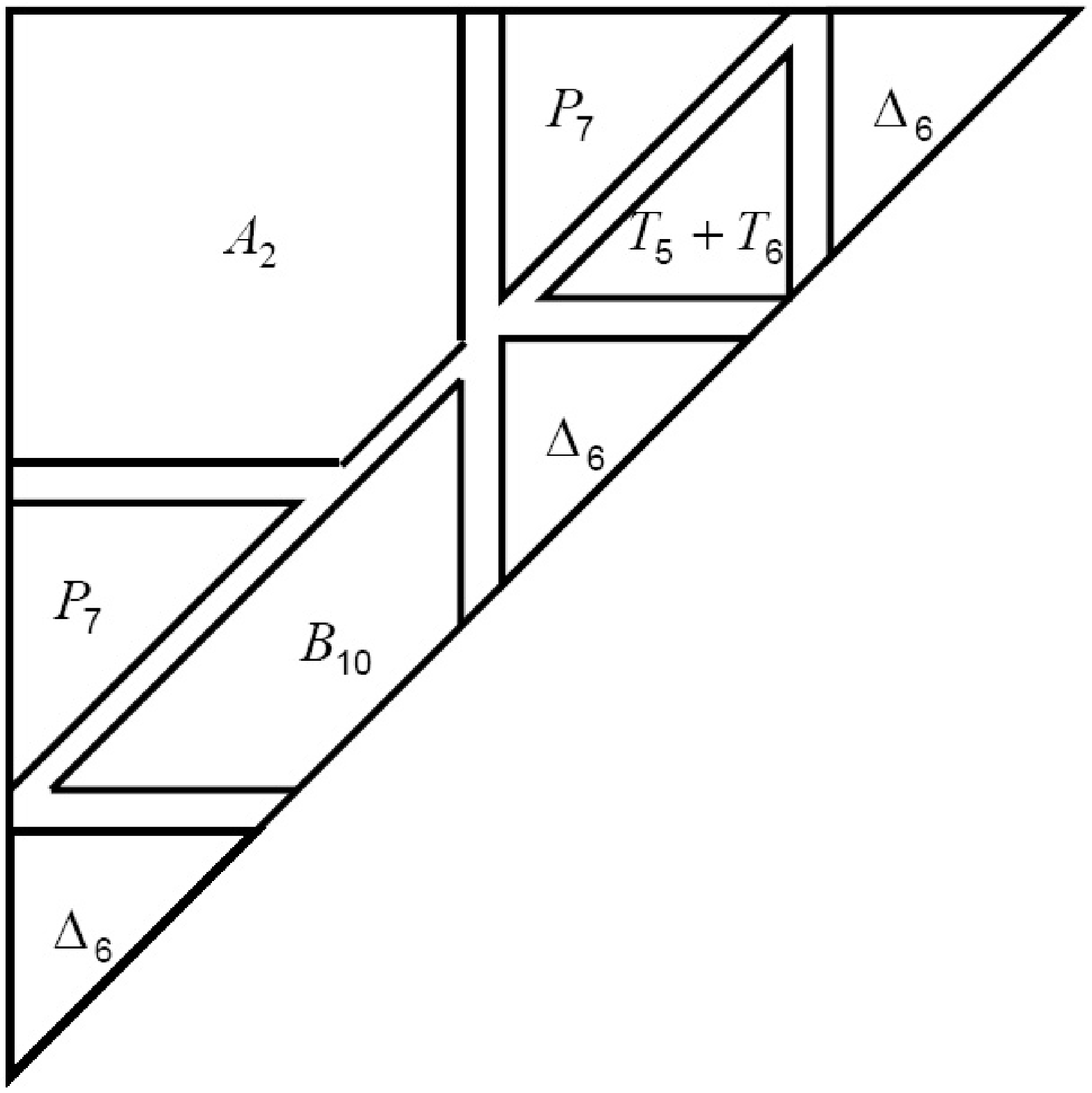}} \qquad\qquad\qquad
\subfigure[Orthogonal projection of $S^7_{34}$.]
{\includegraphics[height=7cm]{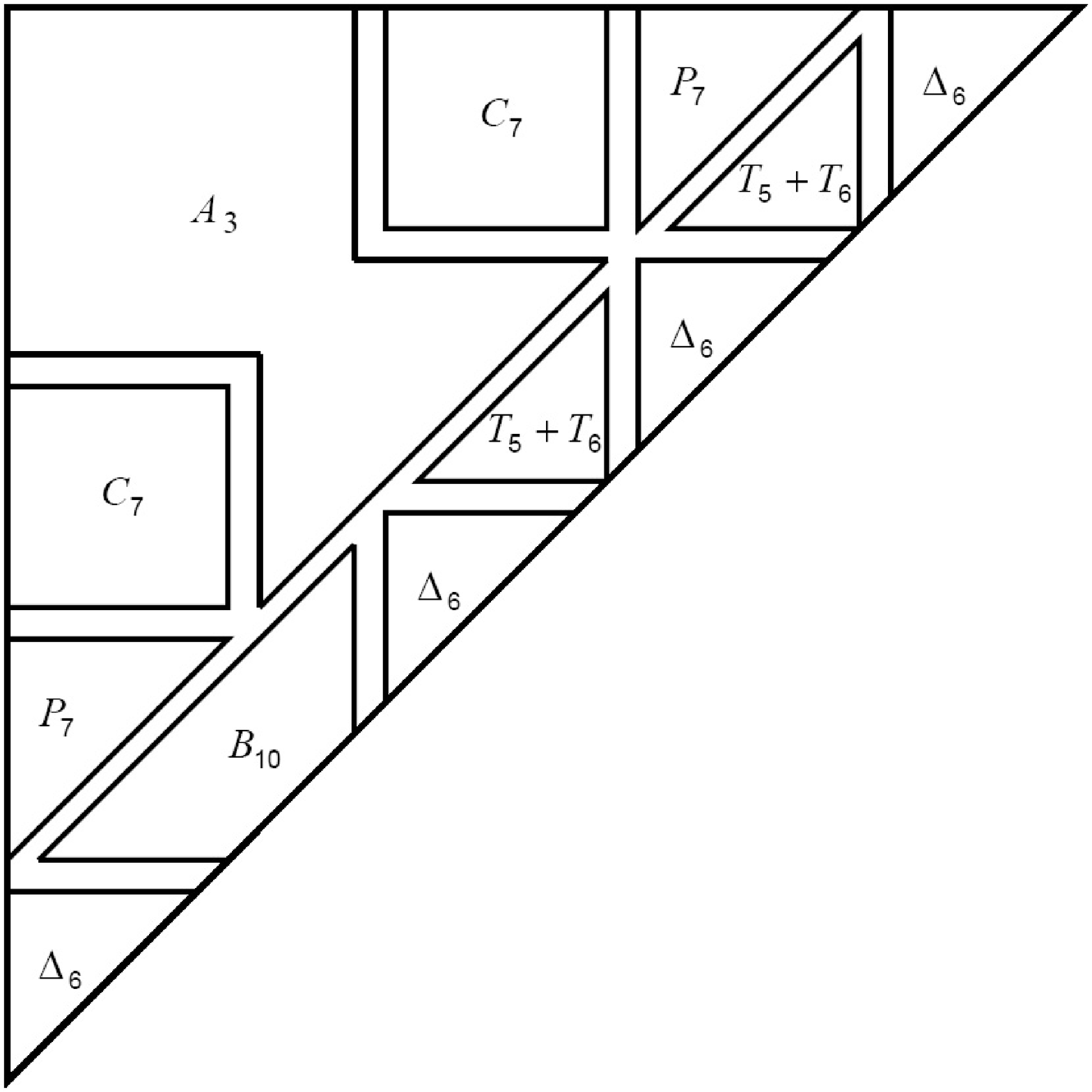}}\caption {\small{Projections
of subdivisions.}}\label{d=26}
\end{figure}
So we get that the decomposition in polytopes is the following (see Figure \ref{d=26}(a) and Figure \ref{d=26}(b)):$$\D_d=\D_{d-8}+A_k+2(k-2)\cdot C_7+2\cdot
P_7+(k+1)\cdot \D_6+B_{10} +(k-1)\cdot (T_5+T_6).$$
The last expression corresponds to the following defectivities:
$$n_{d}=n_{d-8}+(\al_k+1)
+2(k-2)\cdot 128
+2\cdot 72+(k+1)\cdot 21+63+(k-1)\cdot (14+21).$$
Hence we conclude that $V_{3,10+8k}$ is not $n_{10+8k}$--defective.

\subsection{Configuration for d\,=\,12+8k}
\begin{figure}[h]
\centering
\includegraphics[height=4cm]{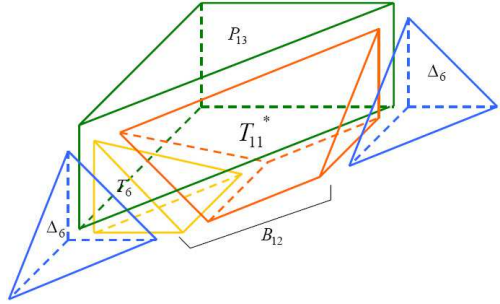}
\caption{\small{Subdivision of $S^7_{20}$.}}\label{d=20}
\end{figure}

\begin{figure}[h] \centering
\subfigure[The cut of $S^5_9$ in $T_{11}^*$.]
{\includegraphics[height=3.5cm]{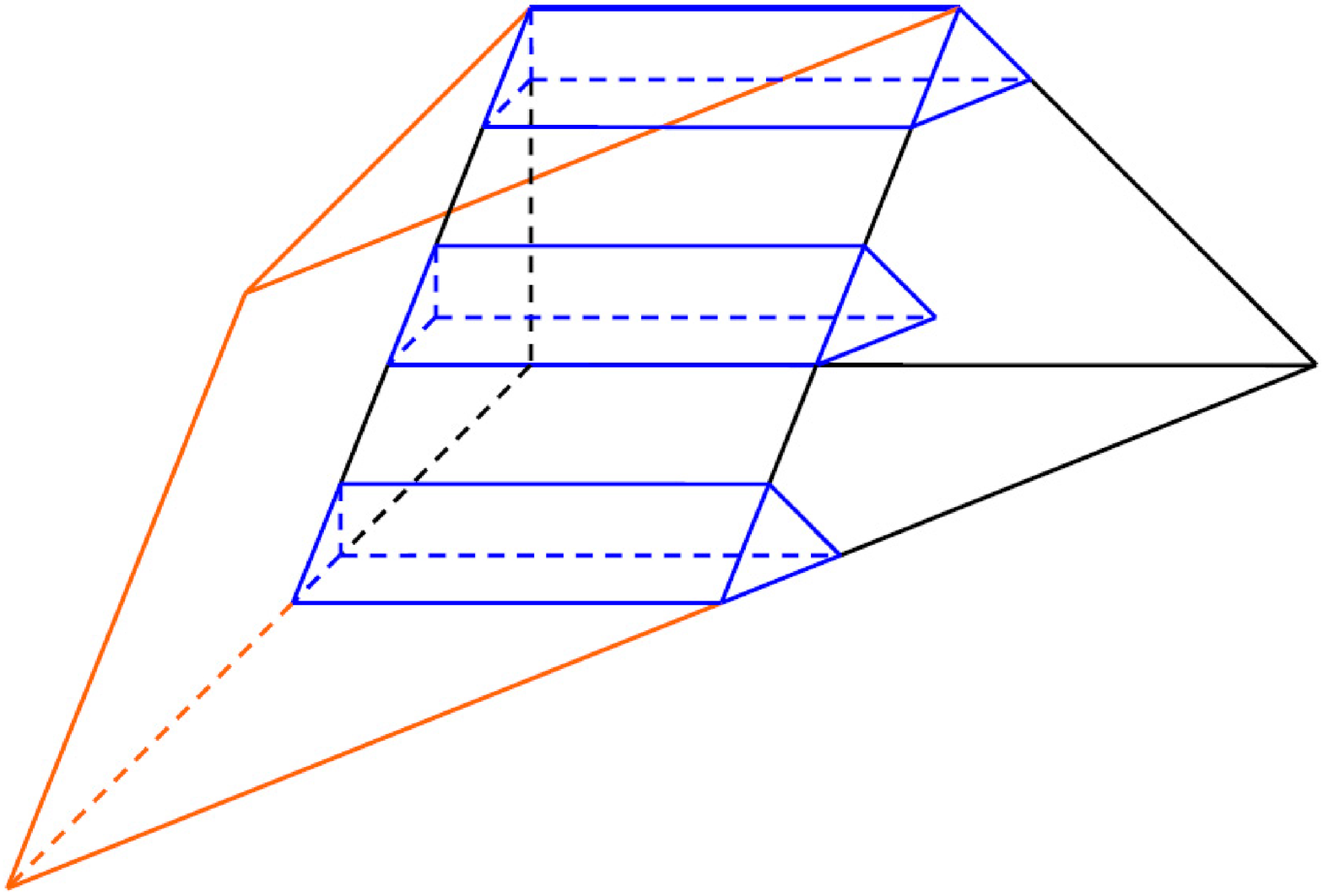}} \qquad\qquad
\subfigure[Orthogonal projection of $A_2$.]
{\includegraphics[height=2.5cm]{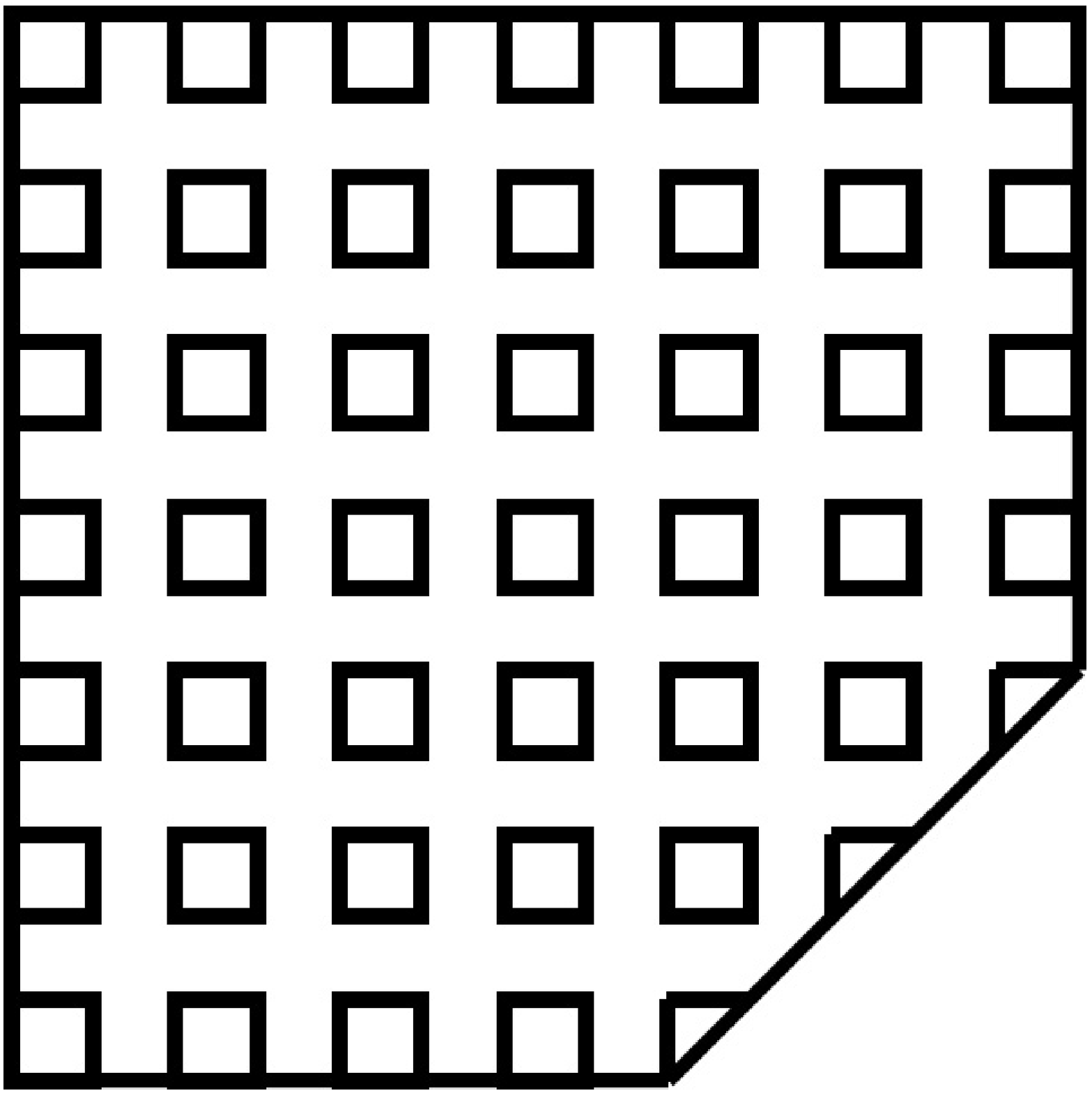}}\qquad\qquad
\subfigure[Orthogonal projection of $A_3$.]
{\includegraphics[height=3.8cm]{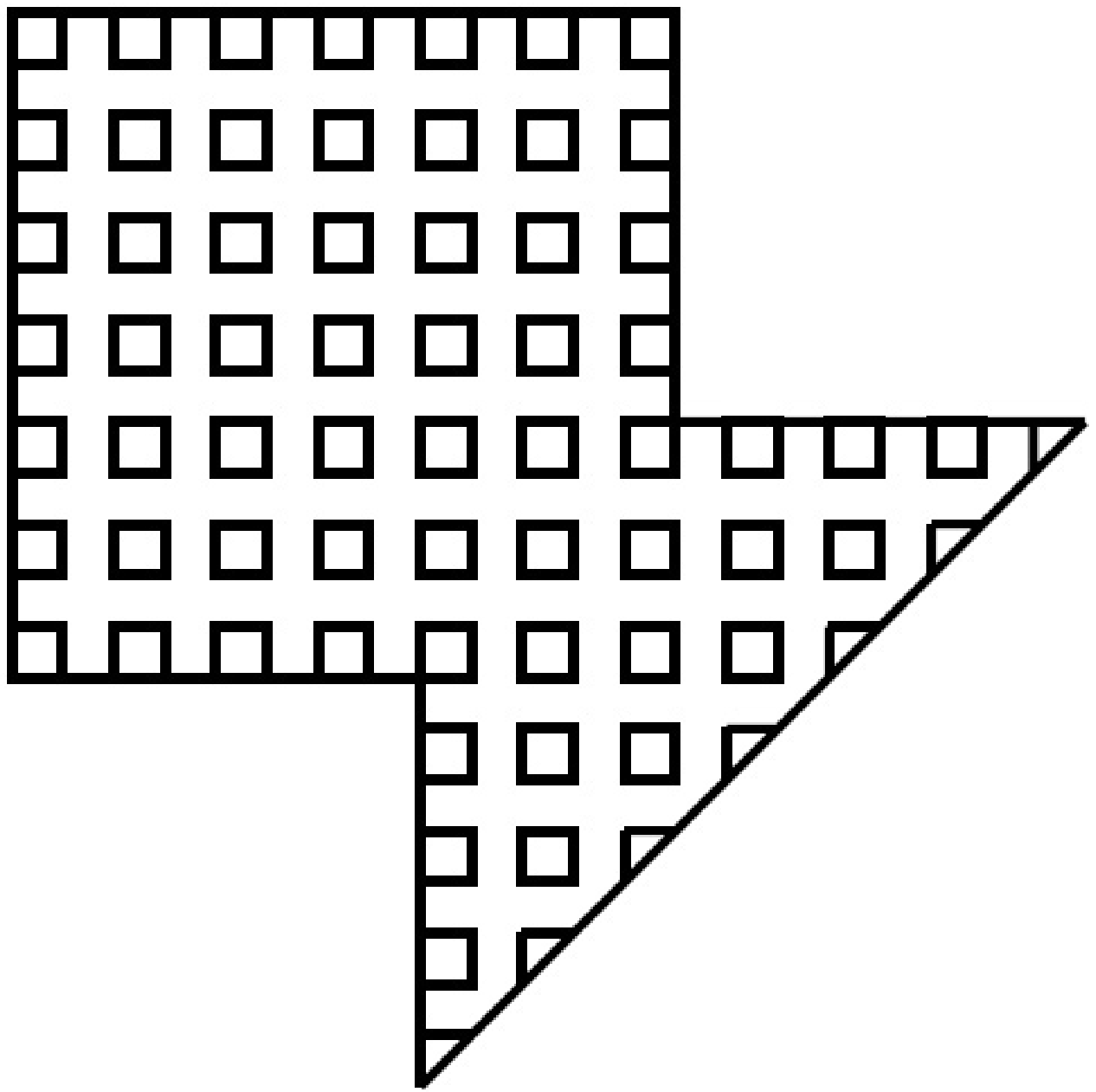}}\caption {\small{}}\label{taglio11}
\end{figure}

\begin{figure}[h]
\centering
\includegraphics[height=7.5cm]{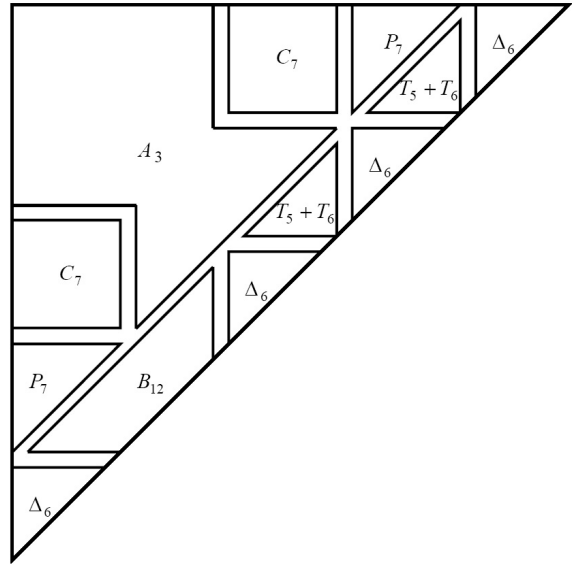}
\caption{\small{Orthogonal projection of $S^7_{36}$.}}\label{d=36}
\end{figure}

This is the last configuration, and it is similar to the one
proposed for $d=10+8k$; again we start examining the case $k=1$
corresponding to $d=20$. We decompose $\D_{20}$ into one tetrahedron
$\D_{12}$ and one $S_{20}^7$. The subdivision of $S_{20}^7$ is the
following: one prism of height $7$ with triangular basis of cathetus
$13$, that we call $P_{13}$; in the front lie $2$ copies of $\D_6$,
one block $B_{12}$ composed of a $T_6$ and a solid $T_{11}^*$ obtained
removing a piece from $\D_{11}$ as in Figure \ref{d=20}.
As we already did for $T^*_7$ and $T^*_9$, we can think of $T^*_{11}$ as obtained by cutting a $\D_{11}$. So we get a solid composed of a $\D_5$ leaning on a layer $S^1_6$ in turn leaning on a $S^5_{11}$ that has been cut as in Figure \ref{taglio11}(a), where we remove the part in black and blue.
\noindent Using again the triangulation shown in
Figure \ref{taglio7}(b) behind the removed blue parts, we see that $T_{11}^*$ is not $55$--defective,
and this implies that the block $B_{12}$
is not $76$--defective (where $21+(91-35)=77$). The decomposition in
polytopes is the following:
$$\D_{20}=\D_{12}+P_{13}+2\cdot \D_6+B_{12}, $$ that is to
say $$n_{20}=n_{12}+210+2\cdot 21+77=442.$$
Then $V_{3,20}$ is not $n_{20}$--defective.
In order to discuss the general case $d=12+8k$, we consider the
subdivision of $\D_d$ in
a $\D_{d-8}$ plus a $S^7_d$. Then we insert in $S_d^7$ a block
$A_k$, whose configuration we will discuss in the following.
For $k=2$, the block $A_2$ contains $43$ disjoint columns of height $7$ and square basis of
side $1$, and $3$ columns $\gamma_7$, hence it is not $361-$defective. The orthogonal projections
of these columns can be viewed in Figure \ref{taglio11}(b). For $k\geq 3$ the block $A_k$ is different: in Figure \ref{taglio11}(c) we see the projection of $A_3$. Looking at this picture, again we see that $A_k$ contains columns of height $7$ and side $1$ (not $7$--defective),
and columns $\gamma_7$ (not $5$--defective). In particular it contains
$\left(46+\sum_{i=4}^{4k-6}i\right)$ columns of the first type, and $(4k-5)$
columns $\gamma_7$. As well as in the previous case we set:
$\al_2+1:=362$, and for $k\geq 3$, $\al_k+1:=8\cdot\left(46+\sum_{i=4}^{4k-6}i\right)+6\cdot (4k-5)$, hence we have that
$A_k$ is not $\al_k$--defective for any $k\geq 2$. In addition to
$A_k$, in $S^7_d$ we insert: $2(k-2)$ copies of $C_7$, $2$ copies of
$P_7$, and in the front stripe there are a block $B_{12}$, $(k+1)$
copies of $\D_6$ plus $(k-1)$ copies of $(T_5+T_6)$. Summarizing, we
have that:
$$\D_d=\D_{d-8}+A_k+2(k-2)\cdot C_7+2\cdot P_7+
(k+1)\cdot \D_6+B_{12}+(k-1)\cdot (T_5+T_6).$$ One example of this
subdivision is available for $S^7_{36}$, whose projection can be
viewed in Figure \ref{d=36}. Substituting in the last expression the
defectivities, we
get:$$n_{d}=n_{d-8}+(\al_k+1)+2(k-2)\cdot
128+2\cdot 72+(k+1)\cdot 21+77+(k-1)\cdot (14+21).$$

\FloatBarrier
\noindent Hence we have that $V_{3,12+8k}$ is not $n_{12+8k}$--defective.
This completes the proof of Alexander-Hirschowitz's theorem in dimension $3$.

\begin{remark}
{\rm Let us observe that in fact we have proved something more: we have decomposed the polytopes $\D_d$ corresponding to $V_{3,d}$, into smaller subpolytopes, and in order to prove the not--defectivity of $V_{3,d}$, we have shown that some of those subpolytopes correspond to not--defective toric varieties; this is further information due to this approach.}
\end{remark}


\begin{thebibliography}{alp}
\bibitem[BO07]{bib:OttBra}M. C. Brambilla, G. Ottaviani, On the Alexander-Hirschowitz Therem, Journal of Pure and Applied Algebra, Vol. 212 (2008), no. 5, 1229-1251

\bibitem[Ch02]{bib:cha}K. A. Chandler, A breef proof of a maximal rank theorem for generic double points in projective space, Trans. Amer. Math. Soc. 353 (2001), no. 5, 1907-1920

\bibitem[Ci01]{bib:interp}C. Ciliberto, Geometric Aspects of
Polynomial Interpolation in More Variables and of Waring's Problem,
Progress in Mathematics, Vol 201, 2001 Birkh\"{a}user Verlag
Basel/Switzerland

\bibitem[CM06]{bib:CC}C. Ciliberto, O. Dumitrescu e R. Miranda,
 Degenerations of the Veronese and Applications, 2006, preprint

\bibitem[DC05]{bib:Cox}D. A. Cox, LECTURES ON TORIC VARIETIES,
2005, indirizzo URL:
http://www.cs.amherst.edu/~dac/lectures/coxcimpa.pdf
\bibitem[JD06]{bib:draisma}Jan Draisma, A tropical approach to secant
dimensions, Mathematics, abstract math.AG/0605345, 2006

\bibitem[EH00]{bib:EisenbudHarris}D. Eisenbud, J. Harris, The
Geometry of Schemes, Springer-Verlag, New York, 2000

\bibitem[Fu93]{bib:Ful}W. Fulton, Introduction to Toric
Varieties, Princeton University Press, Princeton, New Jersey, 1993

\bibitem[GR90]{bib:Gunning}R. C. Gunning, Introduction to
holomorphic functions of several variables, Wadsworth \& Cole
Mathematics Series, Belmont, California, 1990

\bibitem[Ha97]{bib:Hart}R. Hartshorne, Algebraic Geometry,
Springer--Verlag, New York Berlin Heidelberg, 1997

\bibitem[Hu01]{bib:shengda}S. Hu, Semi-Stable Degenerations of
Toric Varieties and Their Hypersurfaces, Mathematics, abstract
math.AG/0110091 v1 8 Oct 2001, University of Wisconsin-Madison, 2001

\bibitem[KY91]{bib:Yang}K. Yang, Complex Algebraic Geometry, An
Introduction to Curves and Surfaces, MARCEL DEKKER, INC., New York,
1991

\bibitem[Mi04]{bib:anacapri}R. Miranda, Anacapri Lectures On Degenerations Of
Surfaces, Anacapri, 2004

\bibitem[Od88]{bib:Oda}T. Oda, Convex Bodies and Algebraic
Geometry--An introduction to the Theory of Toric Varieties,
Springer--Verlag, Berlin Heidelberg, 1988

\bibitem[Sh77]{bib:Shaf} I. Shafarevich, Basic Algebraic Geometry,
Springer--Verlag, Berlin Heidelberg New York, 1977

\bibitem[SS05]{bib:SS} B. Sturmfels, S. Sullivant, Combinatorial secant varietes, math.AC/0506223


\end{thebibliography}
\end{document}